\documentclass[pdftex,10pt,a4paper,oneside]{article}

\usepackage[english]{babel}
\usepackage[utf8]{inputenc}
\usepackage[T1]{fontenc}

\usepackage{mathptmx}

\usepackage{amsmath}
\usepackage{amssymb}
\usepackage{amsthm}

\usepackage{geometry}

\usepackage{newunicodechar}

\newunicodechar{¬}{\lnot}
\newunicodechar{∧}{\land}
\newunicodechar{∨}{\lor}
\newunicodechar{→}{\rightarrow}
\newunicodechar{↔}{\leftrightarrow}
\newunicodechar{∃}{\exists}
\newunicodechar{∀}{\forall}
\newunicodechar{≡}{\equiv}
\newunicodechar{≤}{\leq}
\newunicodechar{≥}{\geq}
\newunicodechar{∈}{\in}
\newunicodechar{∉}{\notin}
\newunicodechar{⊆}{\subseteq}
\newunicodechar{∪}{\cup}
\newunicodechar{∅}{\emptyset}

\usepackage{relsize}

\usepackage{authblk}

\usepackage{titlesec}

\titleformat{\section}[block]{\normalsize\bfseries}{\thesection}{1.0em}{}
\titleformat{\subsection}[runin]{\normalsize\bfseries}{\thesubsection}{0.5em}{}
\titleformat{\paragraph}[runin]{\normalsize\itshape}{\theparagraph}{}{}

\usepackage{tocbibind}

\usepackage{bussproofs}

\renewcommand{\fCenter}{\ensuremath{\Rightarrow}}

\usepackage{bussproofs-extra}

\usepackage[pdftex,unicode,bookmarksnumbered,pagebackref]{hyperref} 

\hypersetup{
 pdftitle={Efficient elimination of Skolem functions in LKʰ},
 pdfauthor={Ján Komara},
 pdfsubject={Mathematical Logic},
 pdfkeywords={
  proof complexity{,} 
  length of proofs{,} 
  Skolem functions{,} 
  eigenvariable condition{,}
  strong locality property{,}
  Henkin constants
 },
}

\hypersetup{
 colorlinks=true,
 allcolors=[rgb]{0,0.5,0.5},
 urlcolor=[rgb]{0.58,0.0,0.83},
 breaklinks=true,
}

\newcommand{\mytexorpdfstring}[2]{\texorpdfstring{#1}{#2}}

\newcommand{\PAR}[1]{\subsection{{#1}\mytexorpdfstring{.}{}}}
\newcommand{\myclaim}{\it}

\numberwithin{equation}{section}

\newcommand{\XEL}[2]{(\ref{eq:#1\if #2\@\else:#2\fi})}

\newcommand{\bMinus}{\stackrel{\textstyle\cdot}{\relbar}}
\newcommand{\aSExp}[2]{\ensuremath{2_{#1}^{#2}}} 

\newcommand{\aSeq}[3][1]{\ensuremath{{#3}_{#1},\dotsc,{#3}_{#2}}}

\let\originalvec\vec
\renewcommand\vec[1]{\originalvec{\kern0pt #1}}
\newcommand*\pvec[1]{\vec{#1}\mkern2mu\vphantom{#1}}

\newcommand{\ARR}[1]{\ensuremath{\Vec{#1}}}   
\newcommand{\ARRr}[1]{\ensuremath{{#1}}}   
\newcommand{\ARRi}[2]{\ensuremath{\Vec{#1}_{#2}}}
   
\newcommand{\ARRpi}[2]{\ensuremath{\pvec{#1}_{#2}'}}

\newcommand{\mEQx}[2]{\ensuremath{\stackrel{\text{#1}}{#2}}}

\newcommand{\nLG}{\ensuremath{{L}}}

\newcommand{\nT}{\ensuremath{\top}}
\newcommand{\nF}{\ensuremath{\bot}}

\newcommand{\nQ}{\ensuremath{Q}}
\newcommand{\nQq}{\ensuremath{Q'}}

\newcommand{\aSubst}[3]{\ensuremath{{#1}_{#2}{[#3]}}}
\newcommand{\aSubstt}[3]{\ensuremath{{#1}_{#2}{\left[#3\right]}}}
\newcommand{\aSubsti}[3]{\ensuremath{{#1}{[#3]}}}

\newcommand{\nHc}{\ensuremath{\text{c}}}
\newcommand{\aHc}[1]{\ensuremath{\nHc_{#1}}}
\newcommand{\nLGC}{\ensuremath{{C}}}
\newcommand{\aLGC}[1]{\ensuremath{{\nLGC_{#1}}}}
\newcommand{\nLGHC}{\ensuremath{\nLG(\nLGC)}}
\newcommand{\nAxH}{\text{H}}
\newcommand{\nAxQ}{\text{Q}}

\newcommand{\aReplace}[2]{\ensuremath{\{ {#1} / {#2} \}}}
\newcommand{\aReplacee}[2]{\ensuremath{\left\{ {#1} / {#2} \right\}}}

\newcommand{\nPK}{\text{PK}}
\newcommand{\nLK}{\text{LK}}
\newcommand{\nLKh}{\ensuremath{\nLK^\text{h}}}
\newcommand{\nGIIIc}{\textbf{G3c}}
\newcommand{\GK}[1]{}
\newcommand{\nLKp}{\ensuremath{\nLK^{+}}}
\newcommand{\nLKpp}{\ensuremath{\nLK^{++}}}
\newcommand{\nD}{\text{D}}

\newcommand{\nL}{\text{L}}
\newcommand{\nR}{\text{R}}
\newcommand{\nAx}{\text{Ax}}
\newcommand{\nLf}{\ensuremath{\nAx{\nF}}}
\newcommand{\nRt}{\ensuremath{\nAx{\nT}}}
\newcommand{\nLn}{\ensuremath{\nL{¬}}}
\newcommand{\nRn}{\ensuremath{\nR{¬}}}
\newcommand{\nLa}{\ensuremath{\nL{∧}}}
\newcommand{\nRa}{\ensuremath{\nR{∧}}}
\newcommand{\nLo}{\ensuremath{\nL{∨}}}
\newcommand{\nRo}{\ensuremath{\nR{∨}}}
\newcommand{\nLi}{\ensuremath{\nL{→}}}
\newcommand{\nRi}{\ensuremath{\nR{→}}}
\newcommand{\nLe}{\ensuremath{\nL{∃}}}
\newcommand{\nRe}{\ensuremath{\nR{∃}}}
\newcommand{\nLu}{\ensuremath{\nL{∀}}}
\newcommand{\nRu}{\ensuremath{\nR{∀}}}
\newcommand{\nCut}{\text{Cut}}

\newcommand{\nRuleAx}{\text{Axiom}}
\newcommand{\nRule}{\text{Rule}}

\newcommand{\nProof}{\ensuremath{\vdash}}
\newcommand{\rProof}{\ensuremath{\mathrel{\nProof}}}
\newcommand{\rProofEq}[1]{\ensuremath{\mathrel{\nProof_{#1}}}}
\newcommand{\rProvableEq}[1]{\ensuremath{\mathrel{\nProof_{#1}}}}
\newcommand{\rProofLe}[1]{\ensuremath{\mathrel{\nProof_{\leq #1}}}}

\newcommand{\rProofLt}[1]{\ensuremath{\mathrel{\nProof_{< #1}}}}

\newcommand{\nP}[1][\pi]{\ensuremath{#1}}
\newcommand{\nPp}[1][\pi]{\ensuremath{#1'}}	
\newcommand{\nPpp}[1][\pi]{\ensuremath{#1''}}	
\newcommand{\aP}[2][\pi]{\ensuremath{{#1}_{#2}}}
\newcommand{\aPp}[2][\pi]{\ensuremath{{#1}_{#2}'}}
\newcommand{\aPpp}[2][\pi]{\ensuremath{{#1}_{#2}''}}

\newcommand{\aDpF}[1]{\ensuremath{\mathrm{d}\left({#1}\right)}}   
\newcommand{\aHrE}[1]{\ensuremath{\mathrm{rk}({#1})}}   
\newcommand{\aHtP}[1]{\ensuremath{\left|{#1}\right|}}   
\newcommand{\aCrP}[1]{\ensuremath{\mathrm{r}\left({#1}\right)}}   

\newcommand{\aWkP}[3]{\ensuremath{{#1}\left[{#2} \fCenter {#3}\right]}}   

\newcommand{\aSeqc}[3][1]{\ensuremath{{#3}_{#1}\dotsc{#3}_{#2}}}
\newcommand{\aSeqq}[4][1]{\ensuremath{{#3}_{#1}^{#4},\dotsc,{#3}_{#2}^{#4}}}
\newcommand{\aSeqnq}[4][1]{\ensuremath{{#3}_{#1},\dotsc,{#3}_{#2}}}
\newcommand{\aSeqnqdr}[5][1]{\ensuremath{{#3}_{#1}{#5},\dotsc,{#3}_{#2}{#5}}}
\newcommand{\aSeqdd}[3][1]{\ensuremath{{#3}_{#1},..\,,{#3}_{#2}}}
\newcommand{\aSeqcdd}[3][1]{\ensuremath{{#3}_{#1}..{#3}_{#2}}}
\newcommand{\aSeqnqdd}[4][1]{\ensuremath{{#3}_{#1},..,{#3}_{#2}}}

\newcommand{\nAwff}{\ensuremath{\Omega_{f}}}
\newcommand{\nAwfe}{\ensuremath{\Omega_∃}}

\newcommand{\nDr}{\ensuremath{\theta}}

\newcommand{\nIq}{\ensuremath{i}}
\newcommand{\aT}[3]{\ensuremath{{#1}_{#2}}}

\newcommand{\nE}{\ensuremath{e_A}}
\newcommand{\aE}[2]{\ensuremath{\nE({#1})}}

\begin{document}


\title{
Efficient elimination of Skolem functions in \nLKh
}

\author{
Ján Komara
\thanks{email: komara@fmph.uniba.sk}
}

\affil{
Department of Applied Informatics,
\\
Faculty of Mathematics, Physics and Informatics,
\\
Comenius University in Bratislava
}

\date{}

\maketitle

\begin{abstract}

Elimination of a single Skolem function in pure logic
increases the length of proofs only linearly.
The result is shown for derivations with cuts
that are free for the Skolem function 
in a sequent calculus with strong locality property.

\end{abstract}



\section{Introduction}
\label{sc:in}

We look at two questions, which are closely related to 
Pudlák's Problem~22~\cite{qq:CloteKrajicek:93}:
\begin{quote}
Given a proof of a formula from the axiom $∀\ARR{x} A[\ARR{x},f(\ARR{x})]$,
where $f$ is a new function symbol,
find a proof of the formula from the axiom $∀\ARR{x} ∃y A[\ARR{x},y]$.
What is the complexity of the elimination of $f$?
\end{quote}

\begin{quote}
Given a proof of a formula $∃\ARR{x} A[\ARR{x},g(\ARR{x})]$,
where $g$ is a new function symbol,
find a proof of the formula $∃\ARR{x} ∀y A[\ARR{x},y]$.
What is the complexity of the elimination of $g$?
\end{quote}
The Skolem functions $f$ and $g$ are called 
a witnessing function for $∀\ARR{x} ∃y A$ and 
a counterexample function to $∃\ARR{x} ∀y A$, respectively.

Skolem functions play an important role both in proof theory and in automatic theorem proving.
Skolemization with witnessing functions is a process that transforms a formula 
(in negation normal form) to an equisatisfiable formula
by replacing all its existential quantifiers with Skolem witnessing functions. 
We are interested here in the reverse process:
given a proof of a formula from the Skolemization of an axiom,
what is the length of the shortest proof of the same formula from the original axiom?
Skolemization with counterexample functions is used 
for transforming formulas to validity-equivalent formulas;
the process, sometimes called Herbrandization, is dual to the previous one.

The complexity of general methods for eliminating Skolem functions from proofs
is at least exponential
(see 
\cite{%
DBLP:journals/jsyml/BaazHW12,%
maehara1955,%
qq:Shoenfield:67,%
qq:Takeuti:87%
}%
),
some of them are even superexponential.
A partial positive solution is given by 
Avigad \cite{DBLP:journals/tocl/Avigad03}: 
theories strong enough to code finite functions can eliminate 
Skolem functions in polynomial time.
But the general problem for predicate calculus either with or without equality
is still an open question.

As an example, consider Maehara's method \cite[Lemma 8.11]{qq:Takeuti:87}
applied to cut-free \nLK-derivations.
A cut-free proof \nP\ of a formula $∀x A[x,f(x)] → B$
is transformed to a proof \nPp\ \emph{with cuts} 
of a formula $∀x ∃y A[x,y] → B$
so that the length of the proof \nPp\ is 
polynomialy bounded by that of the proof \nP.
Cuts are necessary to restore eigenvariable condition.
The subsequent cut elimination leads to superexponential increase in proof length.

As Baaz and Ferm{\"{u}}ller \cite{DBLP:journals/jsyml/BaazHW12} 
has already pointed out, proof complexity of the elimination
does not depend so much on Skolemization per se 
but rather on the nature of the eigenvariable condition.
A step in the right direction is the doublet of sequent calculi \nLKp/\nLKpp\ 
introduced in Aguilera and Baaz \cite{DBLP:journals/jsyml/AguileraB19}.
The authors propose a very liberalized form of eigenvariable condition 
in the form of \emph{locally unsound} but 
\emph{globally correct} strong quantifier rules.
They assert, though without a detailed proof, 
that deskolemization of cut-free \nLKpp-proofs is linear
\cite[Prop. 4.8]{DBLP:journals/jsyml/AguileraB19}.
So, to obtain a better (subexponential)  estimate,
we have to address the conundrum of eigenvariable condition first.
That is done in the next section.

In this paper we consider the problem of Skolemization 
for first-order logic without equality.
We prove that elimination of a single Skolem function 
from for derivations with cuts that are free for the Skolem function 
increases the length of such proofs only linearly.
The result is shown for a variant of sequent calculus, 
called here~\nLKh, with strong locality property. 

The paper is organized as follows.
The motivation behind the calculus \nLKh\ is discussed in Section~\ref{sc:eigen}.
Section~\ref{sc:lg} is completely preparatory. 
In Section~\ref{sc:lk} we give a detailed description of \nLKh.
Section~\ref{sc:ex} contains several examples of \nLKh-derivations.
In Section~\ref{sc:pr} we study both the basic properties of \nLKh-proofs and 
the effect of deep replacement on \nLKh-deductions.
Sections~\ref{sc:cut}~and~\ref{sc:skk} contain the main results:
effective cut elimination for \nLKh-proofs 
via a non-Gentzen-style algorithm without resorting to regularization
and efficient elimination of Skolem functions for a subclass of \nLKh-derivations; 
this includes Problem~22~\cite{qq:CloteKrajicek:93} as well.


\section{How to dispense with eigenvariable condition: motivation behind \mytexorpdfstring{\nLKh}{LKʰ}}
\label{sc:eigen}

Gentzen sequent calculus \nLK, which comes in many flavors,
uses free variables (parameters) in strong quantifier rules.
Free variables behave as uninterpreted constants 
added to a (pure) first-order language, 
they cannot be quantified.
Each free variable can stand for any element of the domain of discourse.
Their meaning, however, can be bound by their use in proofs.

Consider, for instance, a strong quantifier rule \nLe\ of an \nLK-proof 
with the principal formula $∃xA$ as shown on the left:
\def\proofSkipAmount{\vskip 0.8ex plus 0.8ex minus 0.4ex}
\begin{prooftree}
\AxiomC{$\aSubst{A}{x}{a}, \GK{∃xA,} \Gamma \fCenter \Delta$}
\LeftLabel{\nLe}
\RightLabel{,}
\UnaryInfC{$∃xA, \Gamma \fCenter \Delta$}
\DisplayProof\quad
\AxiomC{$\Gamma \fCenter \Delta, \GK{∀xA,} \aSubst{A}{x}{a}$}
\LeftLabel{\nRu}
\RightLabel{.}
\UnaryInfC{$\Gamma \fCenter \Delta, ∀xA$}
\end{prooftree}
\def\proofSkipAmount{\vskip 0.8ex plus 0.8ex minus 0.4ex}
The free variable $a$, which belongs to the \nLe-inference, 
stands for an element of the domain 
for which the minor formula \aSubst{A}{x}{a} is true 
(if there is such element).
For the \nLe-rule to be sound the parameter $a$ must not acquired 
a different meaning earlier during the proof.
The simplest way to ensure this is to 
comply with the following condition:
the proper variable $a$ of the \nLe-inference
occur neither in the principal formula $∃xA$ 
nor in the side formulas $\Gamma, \Delta$ (context) of the rule.
The requirement is called the eigenvariable condition for the \nLe-inference.
A similar requirement holds also for the strong quantifier rule \nRu.

The sequent calculi with the standard eigenvariable condition 
have the following \emph{weak locality property:}
validity of each strong quantifier inference step depends 
not only on the form of its principal and minor formulas, 
but \emph{also} on its context.
Weak locality has some drawbacks though. 
The dependency of correctness for strong quantifier rules on their context 
can make structural analysis of formal proofs 
too technical and overcomplicated.
It may even negatively influence the proof-complexity analysis.

Eigenvariable condition can be somewhat relaxed.
The most notable example is Smullyan's Liberalized Rule \nD\ 
for his analytical tableaux \cite{qq:Smullyan:68} 
that inspired so many followers 
\cite{%
qq:DBLP:conf/tableaux/BaazF95,%
qq:DBLP:conf/kgc/BeckertHS93,%
DBLP:journals/jar/CantoneA07,%
DBLP:journals/jar/HahnleS94%
}.
A modern recent example on the subject is the already mentioned
pair of sequent calculi \nLKp/\nLKpp\ \cite{DBLP:journals/jsyml/AguileraB19}.
The distinguished feature of all these systems is 
that weak locality property is lost:
the validity of a strong quantifier inference step in a proof depends 
not only on its context 
but \emph{also} on the whole history of the derivation, 
or at least a part of it.
This makes structural proof analysis of such a calculus even harder!

\medskip

Instead of making the standard eigenvariable condition more liberal,
we have decided to make it more \emph{rigid}.
The idea is borrowed from 
Komara and Voda \cite{qq:KomaraVoda:95:TABLEAUX95}.
We require that the following two conditions are met:
\begin{itemize}
\item
Strong quantifier rules with the same principal formula $\nQ x A$, 
where $\nQ ∈ \{∃,∀\}$,
have the same unique proper variable \aHc{\nQ x A}.
\item
A rank function \aHrE{a} exists
from free variables onto natural numbers 
such that whenever \aHc{\nQ x A} occurs in a quantifier formula $\nQq y B$ we have
\(
\aHrE{\aHc{\nQ x A}} < \aHrE{\aHc{\nQq y B}}
\).
(Think of number \aHrE{a} as the ordinal number of the level 
in which the parameter $a$ 
occurs for the first time in a cumulative hierarchy of free variables.)
\end{itemize}
One such construction of free variables together with a rank function 
is given in Par.~\ref{par:lg:hc} and 
it is used throughout the whole paper.

In our version of sequent calculus, which we call \nLKh,
strong quantifier rules are thus of the form:
\def\proofSkipAmount{\vskip 0.8ex plus 0.8ex minus 0.4ex}
\begin{prooftree}
\AxiomC{$\aSubst{A}{x}{\aHc{∃xA}}, \GK{∃xA,} \Gamma \fCenter \Delta$}
\LeftLabel{\nLe}
\RightLabel{,}
\UnaryInfC{$∃xA, \Gamma \fCenter \Delta$}
\DisplayProof\quad
\AxiomC{$\Gamma \fCenter \Delta, \GK{∀x A,} \aSubst{A}{x}{\aHc{∀xA}}$}
\LeftLabel{\nRu}
\RightLabel{.}
\UnaryInfC{$\Gamma \fCenter \Delta, ∀xA$}
\end{prooftree}
\def\proofSkipAmount{\vskip 0.8ex plus 0.8ex minus 0.4ex}
We allow for the parameters \aHc{∃xA} and \aHc{∀xA}
to occur in the side formulas $\Gamma, \Delta$;
they may violate the standard eigenvariable condition.
\nLKh-rules have the following \emph{strong locality property:}
the validity of each \nLKh-inference step depends only 
on the form of its principal and minor formulas, 
it does \emph{not} depend on its context. 

The free variables \aHc{∃xA} and \aHc{∀xA},
or rather the \emph{Henkin constants} \aHc{∃xA} and \aHc{∀xA}
as we start to call them from now on,
are called a witness for $∃xA$ and a counterexample to $∀xA$, respectively (see \cite{qq:Henkin:49}).
If $∃xA$ is true, then \aHc{∃xA} denotes an element from the domain which, 
when assigned to the variable $x$, makes $A$ true.
If $∀xA$ is false, then \aHc{∀xA} denotes an element from the domain which, 
when assigned to the variable $x$, makes $A$ false.
This intended meaning of Henkin constants is expressed by 
Henkin witnessing and counterexample axioms.
These are formulas of the form:
\begin{gather*}
∃x A → \aSubst{A}{x}{\aHc{∃x A}}
,
\qquad
\aSubst{A}{x}{\aHc{∀x A}} → ∀x A
.
\end{gather*}
First-order structures that make Henkin axioms true are called Henkin structures.

The sequent calculus \nLKh\ is \emph{locally sound} 
because \nLKh-rules preserve validity in Henkin structures
(compare with Aguilera and Baaz \cite{DBLP:journals/jsyml/AguileraB19}).
\nLKh-provable formulas are therefore valid in every Henkin structure.
The sequent calculus \nLKh\ is also \emph{globally correct} 
since every \nLKh-provable pure formula is logically valid
(see Thm.~\ref{par:lg:henkin:pure}).
(A~formula is pure if it does not contain any parameter.)

We have already examined
in Komara and Voda \cite{qq:KomaraVoda:95:TABLEAUX95}
the feasibility of using Henkin constants in a first-order tableaux calculus. 
The tableaux of \cite{qq:KomaraVoda:95:TABLEAUX95}
are dual to the tableaux of Smullyan \cite{qq:Smullyan:68}
as they demonstrate logical validity rather than logical unsatisfiability.
By direct translation of the arguments \cite{qq:KomaraVoda:95:TABLEAUX95} used for tableaux
we obtain completeness of the \nLKh-calculus for granted 
(see Thm(s).~\ref{par:lk:sc} and~\ref{par:lk:sc:pure}).

At the heart of our results lie a novel syntactic \emph{local} transformation:
\emph{deep replace\-ment} of terms.
By a deep occurrence in an expression we mean, besides ordinary occurrence,
also every deep occurrence in the index $QxA$ of each Henkin constant \aHc{QxA}
occurring in the expression.
The replacement lemma~\ref{par:repl:thm}
characterizes \nLKh-proofs that are invariant under deep replacement of terms.
In particular,
deep replacement turns an invariant strong quantifier rule 
into another valid inference of the same kind.

Both the strong locality property and the deep replacement local transformation 
are specific to \nLKh;
they cannot be readily adapted to a proof calculus
with the standard or a liberalized form of eigenvariable condition.
It is therefore surprising that the calculi \nLKh\ and \nLKpp\ of 
Aguilera and Baaz~\cite{DBLP:journals/jsyml/AguileraB19}
can simulate each other linearly (see Thm.~\ref{par:lkpp:thm}),
as \nLKpp\ does not even have 
the weak locality property besides being locally unsound.
Thus, the proof-theoretic results acquired for one system can be 
applied almost immediately to the other, and vice versa.


\section{First-order logic}
\label{sc:lg}

We start with a quick review of first-order logic.  
We follow the presentation of Komara and Voda \cite{qq:KomaraVoda:95:TABLEAUX95}.
The reader is referred to Barwise \cite{qq:Barwise:77} for more details.

We assume that the reader is familiar with the class of 
elementary recursive functions (see \cite{qq:Rose:82}).
The superexponentiation function \aSExp{n}{x} inductively defined by
$\aSExp{0}{x} = x$ and $\aSExp{n+1}{x} = 2^{\aSExp{n}{x}}$
is an example of non-elementary function.

\PAR{First-order languages}

Logical symbols of first-order languages
include the full set of propositional connectives
$\nT\,\text{(true)}$, $\nF\,\text{(falsehood)}$, $¬$, $∧$, $∨$, $→$, 
and quantifiers $∃$, $∀$;
binary connectives are right associative and 
listed in the order of decreasing precedence.
Each first-order language is fully given 
by the set of its \emph{non-logical symbols:}
constants, function symbols and predicate symbols.

\emph{Semiterms}
are built from variables and constants by application of function symbols. 
Terms are \emph{closed} semiterms; they do not contain variables.
\emph{Semiformulas}
are built from atomic semiformulas 
by application of propositional connectives and quantifiers.
Formulas are \emph{closed} semiformulas; all their variables are bound.
We use letters $r, s, t$ and $A, B, C$ 
to stand for semiterms and semiformulas, respectively.

By expression of a first-order language \nLG\ (\nLG-expression for short) 
we mean either semiterm or semiformula of \nLG.
We use $E_1 ≡ E_2$ as the syntactic identity over expressions of the same kind.

We write \ARR{x} in contexts like $f(\ARR{x})$, where $f$ is an $n$-ary
function symbol, as an abbreviation for a sequence of $n$ variables \aSeq{n}{x},
we have
\(
f(\ARR{x}) ≡ f(\aSeq{n}{x})
\).
Generally, $f(\ARR{t})$ is an abbreviation for $f(\aSeq{n}{t})$,
where \ARR{t} is the sequence \aSeq{n}{t} of terms. 
Semiterms of the form $f(\ARR{t})$ are called \emph{$f$-semiterms.}

By \aSubst{A}{x}{t} we denote the \emph{substitution} of a variable $x$ 
by a term $t$ for every free occurrence of the variable $x$ in a semiformula $A$.
The simultaneous substitution \aSubst{A}{\ARR{x}}{\ARR{t}} is defined analogously. 
We write \aSubsti{A}{\ARR{x}}{\ARR{x}} 
to indicate every free occurrence of (pairwise distinct) variables  \ARR{x} in $A$, and 
\aSubsti{A}{\ARR{x}}{\ARR{t}} as an abbreviation for \aSubst{A}{\ARR{x}}{\ARR{t}}.
A similar notation convention holds also for semiterms.

The notions of \emph{subsemiterm} and \emph{subsemiformula} have a standard definition.
By \emph{subformula} we mean to be that of a subformula in the sense of Gentzen.
For instance, if $∀x B$ or $∃x B$ is a subformula of $A$, 
then so is each (closed) instance \aSubst{B}{x}{t} of $B$.

The \emph{depth} \aDpF{A} of a semiformula $A$ 
is the depth of its tree representation:
\begin{gather*}
\aDpF{A} = 1, \quad \text{$A$ is atomic, \nT\ or \nF},
\\
\aDpF{A ∧ B} = \aDpF{A ∨ B} = \aDpF{A → B} = \max(\aDpF{A},\aDpF{B})+1
\\
\aDpF{¬A} = \aDpF{∃xA} = \aDpF{∀xA} = \aDpF{A}+1.
\end{gather*}

\PAR{Henkin constants and Henkin expansion}
\label{par:lg:hc}

Let \nLG\ be a first-order language. Let further
\(
\aLGC{0}, \aLGC{2}, \aLGC{2}, \ldots, \aLGC{n}, \aLGC{n+1}, \ldots
\)
be an infinite sequence of sets of first-order constants 
defined inductively as follows
\begin{align*}
\aLGC{0} & = \emptyset
\\
\aLGC{n+1} & = 
\aLGC{n} ∪
\left\{ 
\aHc{QxA} \mid 
\text{
$Q ∈ \{∃,∀\}$,
$QxA$ is a formula of $\nLG ∪ \aLGC{n}$,
$\aHc{QxA} \notin \aLGC{n}$
} 
\right\}
\!.
\end{align*}
The set \aLGC{n+1} is thus obtained from the set \aLGC{n} 
by addition of a new constant \aHc{QxA}
for every quantifier formula $QxA$ of the first-order language $\nLG \cup \aLGC{n}$
provided the constant is not already in \aLGC{n}.
The set of all Henkin constants \nLGC\ for \nLG\ is defined as their union:
\[
\nLGC = \bigcup_{n = 0}^\infty \aLGC{n}
.
\]
We refer to the quantifier formula $QxA$ in \aHc{QxA} 
as the \emph{index} of the Henkin constant. 
We also say that the Henkin constant \aHc{QxA} \emph{belongs} to the formula $QxA$.
We use the lowercase letters $a, b, c, \ldots$ to denote Henkin constants.

We obtain the Henkin expansion \nLGHC\ of a first-order language \nLG\ 
by adding all Henkin constants \nLGC\ to the language \nLG:
\[
\nLGHC = \nLG \cup \nLGC
.
\]
\emph{Pure} semiterms and semiformulas of \nLGHC\ are 
expressions in which no Henkin constant occur; 
in fact, they are themselves expressions of the (pure) language \nLG.

The \emph{rank} of a Henkin constant $a$, denoted by \aHrE{a}, is the minimal number $n$ 
such that the constant belongs to the set \aLGC{n}.
It is obvious that
if a Henkin constant \aHc{\nQ x A} occurs in a quantifier \nLGHC-formula $\nQq y B$,
then 
\(
\aHrE{\aHc{\nQ x A}} < \aHrE{\aHc{\nQq y B}}
\).

\PAR{Henkin and quantifier axioms}

Let \nLG\ be a first-order language and \nLGHC\ its Henkin expansion.
Henkin axioms for \nLG\ are \nLGHC-formulas of the form
\begin{gather*}
∃x A → \aSubst{A}{x}{\aHc{∃x A}}
,
\qquad
\aSubst{A}{x}{\aHc{∀x A}} → ∀x A
.
\end{gather*}
The first formula is called \emph{witnessing axiom} and
the second \emph{counterexample axiom}.
The Henkin constants \aHc{∃x A} and \aHc{∀x A} are called 
a \emph{witness} for $∃x A$ and a \emph{counterexample} to $∀x A$, respectively.

Quantifier axioms for \nLG\ are \nLGHC-formulas of the form 
($t$ is arbitrary \nLGHC-term)
\begin{gather*}
\aSubst{A}{x}{t} → ∃x A
,
\qquad 
∀x A → \aSubst{A}{x}{t}
.
\end{gather*}
It is obvious that quantifier axioms are logically valid formulas.

By a \emph{Henkin structure} we mean every first-order structure for \nLGHC\
that is a model for the set consisting of all Henkin axioms for \nLG.

\PAR{The reduction to propositional logic}

The combination of the next two theorems reduces the problem of 
recognizing logical validity for pure formulas
to conceptually much simpler problem of 
recognizing a certain (propositional) tautology.
From this one gets without difficulty 
soundness and completeness for ordinary first-order proof calculi.
Smullyan \cite{qq:Smullyan:68} calls the idea central to predicate logic
in the form of his Fundamental Theorem.
A modern and very readable presentation is by Barwise \cite{qq:Barwise:77}.

\PAR{Theorem  \textnormal{(Henkin reduction)}}
\label{par:lg:henkin:reduction}
{\myclaim

A formula $A$ (not necessarily pure) is valid in every Henkin structure 
if and only if 
there are finitely many Henkin and quantifier axioms $A_1, \ldots, A_n$ 
such that the implication
\(
A_1 ∧ \cdots ∧ A_n → A
\)
is a tautology.

}

\begin{proof}
It follows from the proof of Main Lemma 4.8 and 
Compactness Theorem for Propositional Logic 4.2 
in \cite{qq:Barwise:77}.
\end{proof}

\PAR{Theorem \textnormal{(Pure formula reduction)}}
\label{par:lg:henkin:pure}
{\myclaim

A pure formula is logically valid  if and only if 
it is valid in every Henkin structure.

}

\begin{proof}
By adapting the proof of Lemma 4.7 in \cite{qq:Barwise:77} 
to include also the Henkin counterexample constants.
\end{proof}


\section{Sequent calculus \mytexorpdfstring{\nLKh}{LKʰ}}
\label{sc:lk}

The axiom system \nLKh\ is a variant of the sequent calculus \nGIIIc.
The distinguished feature of \nGIIIc\ is 
that sequents are pairs of multisets of formulas.
Weak structural rules such as contraction and weakening
are thus absorbed into logical inference  rules and axioms
(see \cite{qq:TroelstraSchwichtenberg:PT:00} for details).
The main difference between \nGIIIc\ and our method of sequent calculus is that 
\nLKh\ uses Henkin constants as free variables.

Throughout the section we assume that
a first-order language \nLG\ and its Henkin expansion \nLGHC\ is fixed.
All formulas and terms, unless otherwise specified, 
are first-order expressions of \nLGHC.

\PAR{Finite multisets}

By finite multisets we mean 
finite unordered collections of elements with repetitions. 
So a finite multiset is like an ordinary set 
only it may contain some elements with multiple occurrences.
We use Greek capitals
$\Gamma, \Delta, \Pi, \Lambda, \Omega$ 
to stand for finite multisets.

For finite multisets we adopt the following notation conventions.
If $A$ is an element, then by the same symbol $A$ 
we denote the finite multiset containing $A$ as its only element.
By $\Gamma, \Delta$ we denote the union of finite multisets $\Gamma$ and $\Delta$.
Multisets $A, \Gamma$ and $\Gamma, A$ are thus 
the result of adding the element $A$ to the finite multiset $\Gamma$.

\PAR{Sequents}

By a sequent $\Gamma \fCenter \Delta$ 
we mean the ordered pair $(\Gamma, \Delta)$ 
of finite (possibly empty) multisets of formulas, 
The sequent arrow $\fCenter$ separates 
the \emph{antecedent} $\Gamma$ of the sequent from its \emph{succedent} $\Delta$.
Together they are called \emph{cedents}.
We use the uppercase letter $S$, possibly subscripted,
as a syntactic variable ranging over sequents.

The meaning of a sequent of the form
$A_1, \ldots, A_m \fCenter B_1, \ldots, B_n$
is equivalent in the meaning to the formula
$A_1 ∧ \cdots ∧ A_m → B_1 ∨ \cdots ∨ B_n$.
An empty conjunction ($m = 0$) is defined to be \nT\ and
an empty disjunction ($n = 0$) is defined to be \nF.
The sequent $\fCenter A$ has thus the same meaning as the formula $A$ and 
the \emph{empty} sequent $\fCenter$ is logically unsatisfiable.

\PAR{Inference rules}

The rules of inference of \nLKh\ are of the form
\def\proofSkipAmount{\vskip 0.8ex plus 0.8ex minus 0.4ex}
\begin{prooftree}
\AxiomC{\phantom{$S_1$}}
\RightLabel{,}
\UnaryInfC{$S$} 
\DisplayProof
\qquad
\AxiomC{$S_1$}
\RightLabel{,}
\UnaryInfC{$S$} 
\DisplayProof
\qquad
\AxiomC{$S_1$}
\AxiomC{$S_2$}
\RightLabel{,}
\BinaryInfC{$S$}
\end{prooftree}
\def\proofSkipAmount{\vskip 0.8ex plus 0.8ex minus 0.4ex}
where $S, S_1, S_2$ are sequents.
\emph{Premises} (\emph{hypotheses}) of a rule are its upper sequents,
the \emph{conclusion} of a rule is its lower sequent.
Nullary rules, that is rules with no premises, are called \emph{axioms}.
The list of valid rules of inference of \nLKh\ follows:

\begin{list}{}{\setlength{\leftmargin}{0in}\setlength{\rightmargin}{0in}}

\item \emph{Axioms}

\def\proofSkipAmount{\vskip 0.8ex plus 0.8ex minus 0.4ex}
\begin{prooftree}
\AxiomC{}
\LeftLabel{\nAx}
\RightLabel{($A$ is atomic),}
\UnaryInfC{$A, \Gamma \fCenter \Delta, A$} 
\DisplayProof
\quad
\AxiomC{}
\LeftLabel{\nLf}
\RightLabel{,}
\UnaryInfC{$\nF, \Gamma \fCenter \Delta$}
\DisplayProof
\quad
\AxiomC{}
\LeftLabel{\nRt}
\RightLabel{.}
\UnaryInfC{$\Gamma \fCenter \Delta, \nT$}
\end{prooftree}
\def\proofSkipAmount{\vskip 0.8ex plus 0.8ex minus 0.4ex}

\item \emph{Propositional rules}

\def\proofSkipAmount{\vskip 0.8ex plus 0.8ex minus 0.4ex}
\begin{prooftree}
\AxiomC{$\GK{¬A,} \Gamma \fCenter \Delta, A$}
\LeftLabel{\nLn}
\RightLabel{,}
\UnaryInfC{$¬A, \Gamma \fCenter \Delta$}
\DisplayProof\quad
\AxiomC{$A, \Gamma \fCenter \Delta\GK{, ¬A}$}
\LeftLabel{\nRn}
\RightLabel{,}
\UnaryInfC{$\Gamma \fCenter \Delta, ¬A$}
\end{prooftree}
\def\proofSkipAmount{\vskip 0.8ex plus 0.8ex minus 0.4ex}

\def\proofSkipAmount{\vskip 0.8ex plus 0.8ex minus 0.4ex}
\begin{prooftree}
\AxiomC{$A, B, \GK{A ∧ B,} \Gamma \fCenter \Delta$}
\LeftLabel{\nLa}
\RightLabel{,}
\UnaryInfC{$A ∧ B, \Gamma \fCenter \Delta$}
\DisplayProof\quad
\AxiomC{$\Gamma \fCenter \Delta, \GK{A ∧ B,} A$}
\AxiomC{$\Gamma \fCenter \Delta, \GK{A ∧ B,} B$}
\LeftLabel{\nRa}
\RightLabel{,}
\BinaryInfC{$\Gamma \fCenter \Delta, A ∧ B$}
\end{prooftree}
\def\proofSkipAmount{\vskip 0.8ex plus 0.8ex minus 0.4ex}

\def\proofSkipAmount{\vskip 0.8ex plus 0.8ex minus 0.4ex}
\begin{prooftree}
\AxiomC{$A, \GK{A ∨ B,} \Gamma \fCenter \Delta$}
\AxiomC{$B, \GK{A ∨ B,} \Gamma \fCenter \Delta$}
\LeftLabel{\nLo}
\RightLabel{,}
\BinaryInfC{$A ∨ B, \Gamma \fCenter \Delta$}
\DisplayProof\quad
\AxiomC{$\Gamma \fCenter \Delta, \GK{A ∨ B,} A, B$}
\LeftLabel{\nRo}
\RightLabel{,}
\UnaryInfC{$\Gamma \fCenter \Delta, A ∨ B$}
\end{prooftree}
\def\proofSkipAmount{\vskip 0.8ex plus 0.8ex minus 0.4ex}

\def\proofSkipAmount{\vskip 0.8ex plus 0.8ex minus 0.4ex}
\begin{prooftree}
\AxiomC{$\GK{A → B,} \Gamma \fCenter \Delta, A$}
\AxiomC{$B, \GK{A → B,} \Gamma \fCenter \Delta$}
\LeftLabel{\nLi}
\RightLabel{,}
\BinaryInfC{$A → B, \Gamma \fCenter \Delta$}
\DisplayProof\quad
\AxiomC{$A, \Gamma \fCenter \Delta, \GK{A → B,} B$}
\LeftLabel{\nRi}
\RightLabel{.}
\UnaryInfC{$\Gamma \fCenter \Delta, A → B$}
\end{prooftree}
\def\proofSkipAmount{\vskip 0.8ex plus 0.8ex minus 0.4ex}

\item \emph{Quantifier rules}

\def\proofSkipAmount{\vskip 0.8ex plus 0.8ex minus 0.4ex}
\begin{prooftree}
\AxiomC{$\aSubst{A}{x}{\aHc{∃x A}}, \GK{∃x A,} \Gamma \fCenter \Delta$}
\LeftLabel{\nLe}
\RightLabel{,}
\UnaryInfC{$∃x A, \Gamma \fCenter \Delta$}
\DisplayProof\quad
\AxiomC{$\Gamma \fCenter \Delta, ∃x A, \aSubst{A}{x}{t}$}
\LeftLabel{\nRe}
\RightLabel{,}
\UnaryInfC{$\Gamma \fCenter \Delta, ∃x A$}
\end{prooftree}
\def\proofSkipAmount{\vskip 0.8ex plus 0.8ex minus 0.4ex}

\def\proofSkipAmount{\vskip 0.8ex plus 0.8ex minus 0.4ex}
\begin{prooftree}
\AxiomC{$\aSubst{A}{x}{t}, ∀x A, \Gamma \fCenter \Delta$}
\LeftLabel{\nLu}
\RightLabel{,}
\UnaryInfC{$∀x A, \Gamma \fCenter \Delta$}
\DisplayProof\quad
\AxiomC{$\Gamma \fCenter \Delta, \GK{∀x A,} \aSubst{A}{x}{\aHc{∀x A}}$}
\LeftLabel{\nRu}
\RightLabel{.}
\UnaryInfC{$\Gamma \fCenter \Delta, ∀x A$}
\end{prooftree}
\def\proofSkipAmount{\vskip 0.8ex plus 0.8ex minus 0.4ex}

\item \emph{Structural rules}

\def\proofSkipAmount{\vskip 0.8ex plus 0.8ex minus 0.4ex}
\begin{prooftree}
\AxiomC{$\Gamma \fCenter \Delta, A$}
\AxiomC{$A, \Gamma \fCenter \Delta$}
\LeftLabel{\nCut}
\RightLabel{.}
\BinaryInfC{$\Gamma \fCenter \Delta$}
\end{prooftree}
\def\proofSkipAmount{\vskip 0.8ex plus 0.8ex minus 0.4ex}

\end{list}
The quantifier rules \nLe\ and \nRu\ are called \emph{strong} quantifier rules, 
the other two quantifier rules are called \emph{weak}.
We say that the Henkin constants \aHc{∃x A} and \aHc{∀x A} belong 
to the strong quantifier rules \nLe\ and \nRu, respectively.

For every rule, the formula introduced into its conclusion 
is called the \emph{principal (main) formula},
formulas from which the principal formula is derived
are called the \emph{minor (auxiliary) formulas},
and all the remaining formulas are called the \emph{side formulas} (\emph{context}).
In the axiom \nAx\ and the cut rule 
both occurrences of $A$ are principal,
in the axioms \nLf\ and \nRt\ 
both occurrences of the propositional constants are principal.
The principal formula of a cut rule is called 
the \emph{cut-formula} of the inference.

\PAR{Proofs}

A proof in \nLKh\ is a rooted labeled finite tree with sequents as its nodes.
The root of the tree, written at the bottom, 
is called \emph{endsequent} and it is the sequent to be proved. 
The remaining nodes of the tree are built by inference rules.
The leaves, at the top of the tree, 
are \emph{initial sequents} inferred by axiom rules.
The inner nodes of the tree are inferred by the remaining inferences.
A proof is called \emph{cut-free} if it does not contain any cut rule.

We write
\(
\nP \rProof \Gamma \fCenter \Delta
\)
if \nP\ is an \nLKh-proof of $\Gamma \fCenter \Delta$,
and
\(
\rProof \Gamma \fCenter \Delta
\)
if there is such a proof.
By an \nLKh-proof of a formula $A$ we mean 
any \nLKh-proof of its corresponding sequent $\fCenter A$.
We say that a sequent
$\Gamma \fCenter \Delta$ is \nLKh-provable from a set of formulas $T$ 
if 
for some subset $\Pi$ of $T$ the sequent
$\Pi, \Gamma \fCenter \Delta$ is \nLKh-provable.

The \emph{length} of a proof \nP\ is the number of sequents in \nP.
We write
\(
\nP \rProofEq{l} \Gamma \fCenter \Delta
\)
if the length of the proof \nP\ is $l$, and
\(
\rProofEq{l} \Gamma \fCenter \Delta
\)
if there is such a proof.
We write
\(
\nP \rProofLe{l} \Gamma \fCenter \Delta
\)
if the length of the proof \nP\ is $\leq l$,
and
\(
\rProofLe{l} \Gamma \fCenter \Delta
\)
if there is such a proof.
The notation
\(
\nP \rProofLt{l} \Gamma \fCenter \Delta
\)
and 
\(
\rProofLt{l} \Gamma \fCenter \Delta
\)
is defined analogously.

The \emph{height} of a proof \nP, denoted by \aHtP{\nP}, is the height of its tree representation.
It is the number of nodes along a maximal branch of the tree.
This includes leaves inferred by axioms.
We thus always have $\aHtP{\nP} > 0$.
Note also that
\(
\aHtP{\nP} ≤ l < 2^{\aHtP{\nP}}
\)
if the length of the proof \nP\ is $l$.

The \emph{rank} of a cut rule is the depth of its cut-formula.
The \emph{cut-rank} (\emph{rank} for short) of a proof \nP, denoted by \aCrP{\nP},
is the supremum of the ranks of all its cut inferences:
\[
\aCrP{\nP} = \sup \{ \aDpF{A} \mid \text{$A$ is a cut-formula of \nP} \}
.
\]
We have $\aCrP{\nP} = 0$ if and only if the proof \nP\ is cut-free.

\PAR{The characterization problem}

We give here a brief outline of soundness and completeness for the \nLKh-calculus.
It is an adaptation of the proof of a similar theorem 
in Komara and Voda \cite{qq:KomaraVoda:95:TABLEAUX95}
for a tableaux-based system with Henkin constants in place of proper variables. 
By \nPK\ we mean here a propositional sequent calculus 
containing all inferences of \nLKh\ except for quantifier rules.

\PAR{Soundness and completeness theorem I}
\label{par:lk:sc}
{\myclaim

A sequent (not necessarily pure) is \nLKh-provable if and only if
it is valid in every Henkin structure.

}

\begin{proof}
Let \nAxH\ and \nAxQ\ be
the set of all Henkin and quantifier axioms, respectively.
A~sequent $\Gamma \fCenter \Delta$ is \nLKh-provable iff,
by the elimination and introduction of quantifier rules 
(see \cite[Thm. 12]{qq:KomaraVoda:95:TABLEAUX95}),
the sequent $\Gamma \fCenter \Delta$ is \nPK-provable 
in the set of axioms $\nAxH \cup \nAxQ$ iff,
by the deduction theorem for propositional logic 
(see \cite[Thm. 10]{qq:KomaraVoda:95:TABLEAUX95}),
for some $\{A_1, \ldots, A_n\} \subseteq \nAxH \cup \nAxQ$
the sequent $A_1, \ldots, A_n, \Gamma \fCenter \Delta$ is \nPK-provable iff,
by the soundness and completeness theorem for \nPK\
(see \cite[Thm. 7]{qq:KomaraVoda:95:TABLEAUX95}),
for some $\{A_1, \ldots, A_n\} \subseteq \nAxH \cup \nAxQ$
the sequent $A_1, \ldots, A_n, \Gamma \fCenter \Delta$ is a tautology iff,
by the Henkin reduction
(see Thm.~\ref{par:lg:henkin:reduction}),
the sequent $\Gamma \fCenter \Delta$ is valid in every Henkin structure.
\end{proof}

\PAR{Soundness and completeness theorem II}
\label{par:lk:sc:pure}
{\myclaim

A pure sequent is \nLKh-provable if and only if
it is logically valid.

}

\begin{proof}
It follows directly from Thm.~\ref{par:lk:sc} and Thm.~\ref{par:lg:henkin:pure}.
\end{proof}

\PAR{The sequent calculus \mytexorpdfstring{\nLKpp}{LK⁺⁺}}
\label{par:lkpp:def}

The proof system \nLKpp\ of
Aguilera and Baaz \cite{DBLP:journals/jsyml/AguileraB19}
is defined like an ordinary \nLK-calculus, 
except that the proper variables of strong quantifier rules 
in each valid \nLKpp-proof \nP\ have the following three properties:
\begin{itemize}
\item
Substitutability:
the proper variables of \nP\ do not appear in the endsequent of \nP.
\item
Side-variable condition:
the relation $<_{\nP}$ defined on the proper variables of \nP\ is acyclic.
The property $a <_{\nP} b$ holds if the principal formula 
of a strong quantifier rule in \nP\
to which the parameter $a$ belongs 
contains the parameter $b$.
\item
Very weak regularity:
whenever a proper variable of \nP\ is also 
the proper variable of another inference,
then both rules have the same principal formula.
\end{itemize}
The next proof reveals that the relationship between \nLKpp\ and \nLKh\ is quite intimate:
\nLKpp-\ and \nLKh-deductions are the same, up to the renaming of proper variables.

\PAR{Theorem \textnormal{(Linear simulation)}}
\label{par:lkpp:thm}
{\myclaim

A pure formula has an \nLKpp-proof with length $l$ if and only if 
it has an \nLKh-proof with length $l$. 

}

\begin{proof}
Each \nLKh-proof \nP\ with pure endsequent is already a valid \nLKpp-derivation 
since the property
\(
\aHc{\nQ x A} <_{\nP} \aHc{\nQq y B}
\)
implies
\(
\aHrE{\aHc{\nQ x A}} > \aHrE{\aHc{\nQq y B}}
\).
Vice versa, if \nP\ is an \nLKpp-proof, 
then we may suppose without loss of generality
that its strong quantifier rules with the same principal formula 
have the same proper variable.
Let $<_{\nP}^{'}$ be the partial order induced by the acyclic relation $<_{\nP}$.
We take a $<_{\nP}^{'}$-maximal parameter $a$ 
that belongs to a quantifier formula $\nQ x A$
and replace it in \nP\ by the Henkin constant \aHc{\nQ x A}.
The process is repeated until we get a valid \nLKh-proof 
with the same endsequent and length as the \nLKpp-proof \nP.
\end{proof}


\section{Examples of \mytexorpdfstring{\nLKh}{LKʰ}-proofs}
\label{sc:ex}

We now give some examples of \nLKh-derivations 
so that the reader gets familiar 
with the use of Henkin constants in formal proofs.
The reader interested primarily in the main result of the paper 
may skip and go directly to the next section.

\PAR{Example}
\label{par:ex:ex1}

The following tree is an \nLKh-proof 
\def\proofSkipAmount{\vskip 0.8ex plus 0.8ex minus 0.4ex}
\begin{prooftree}
\AxiomC{}
\LeftLabel{\nAx}
\UnaryInf$P(a,b), ∀x P(x,b) \fCenter ∃y P(a,y), P(a,b)$
\LeftLabel{\nRe}
\UnaryInf$P(a,b), ∀x P(x,b) \fCenter ∃y P(a,y)$
\LeftLabel{\nLu}
\UnaryInf$∀x P(x,b) \fCenter ∃y P(a,y)$
\LeftLabel{\nLe}
\RightLabel{$\left(b ≡ \aHc{∃y ∀x P(x,y)}\right)$}
\UnaryInf$∃y ∀x P(x,y) \fCenter ∃y P(a,y)$
\LeftLabel{\nRu}
\RightLabel{$\left(a ≡ \aHc{∀x∃y P(x,y)}\right)$}
\UnaryInf$∃y ∀x P(x,y) \fCenter ∀x ∃y P(x,y)$
\LeftLabel{\nRi}
\UnaryInf$\fCenter ∃y ∀x P(x,y) → ∀x ∃y P(x,y)$
\end{prooftree}
\def\proofSkipAmount{\vskip 0.8ex plus 0.8ex minus 0.4ex}
of the formula $∃y ∀x P(x,y) → ∀x ∃y P(x,y)$.
The defining abbreviations for the Henkin constants $a$ and $b$ are shown 
on the right next to the rule they belong.
The inferences in question are the following strong quantifier rules
shown in the full form:
\def\proofSkipAmount{\vskip 0.8ex plus 0.8ex minus 0.4ex}
\begin{prooftree}
\Axiom$∃y ∀x P(x,y) \fCenter ∃y P\!\left(\aHc{∀x∃y P(x,y)},y\right)$
\LeftLabel{\nRu}
\UnaryInf$∃y ∀x P(x,y) \fCenter ∀x ∃y P(x,y)$
\end{prooftree}
\begin{prooftree}
\Axiom$∀x P\!\left(x,\aHc{∃y ∀x P(x,y)}\right) \fCenter ∃y P\!\left(\aHc{∀x∃y P(x,y)},y\right)$
\LeftLabel{\nLe}
\RightLabel{.}
\UnaryInf$∃y ∀x P(x,y) \fCenter ∃y P\!\left(\aHc{∀x∃y P(x,y)},y\right)$
\end{prooftree}
\def\proofSkipAmount{\vskip 0.8ex plus 0.8ex minus 0.4ex}
When the tree is read as a standard \nLK-proof,
then both parameters $a$ and $b$ satisfy the standard eigenvariable condition.

\PAR{Example \textnormal{(Smullyan's drinker paradox)}}
\label{par:ex:ex2}

Consider now an \nLKh-proof 
\def\proofSkipAmount{\vskip 0.8ex plus 0.8ex minus 0.4ex}
\begin{prooftree}
\AxiomC{}
\LeftLabel{\nAx}
\UnaryInf$D(a) \fCenter D(a), ∃x (D(x) → ∀y D(y))$
\LeftLabel{\nRu}
\RightLabel{$\left(a ≡ \aHc{∀yD(y)}\right)$}
\UnaryInf$D(a) \fCenter ∀y D(y), ∃x (D(x) → ∀y D(y))$
\LeftLabel{\nRi}
\UnaryInf$\fCenter D(a) → ∀y D(y), ∃x (D(x) → ∀y D(y))$
\LeftLabel{\nRe}
\UnaryInf$\fCenter ∃x (D(x) → ∀y D(y))$
\end{prooftree}
\def\proofSkipAmount{\vskip 0.8ex plus 0.8ex minus 0.4ex}
of the formula $∃x (D(x) → ∀y D(y))$.
The Henkin constant $a$ belongs to the strong quantifier inference
(shown in unabbreviated form):
\def\proofSkipAmount{\vskip 0.8ex plus 0.8ex minus 0.4ex}
\begin{prooftree}
\Axiom$D\!\left(\aHc{∀y D(y)}\right) \fCenter D\!\left(\aHc{∀y D(y)}\right), ∃x (D(x) → ∀y D(y))$
\LeftLabel{\nRu}
\RightLabel{.}
\UnaryInf$D\!\left(\aHc{∀y D(y)}\right) \fCenter ∀y D(y), ∃x (D(x) → ∀y D(y))$
\end{prooftree}
\def\proofSkipAmount{\vskip 0.8ex plus 0.8ex minus 0.4ex}
When the tree is read as a standard \nLK-proof,
then the proper variable $a$ of the strong quantifier rule \nRu\
violates the standard eigenvariable condition.
This is because $a$ occurs in the side formula $D(a)$ 
of the conclusion of the \nRu-inference.
This kind of violation is still sound because 
it satisfies the prerequisites of Smullyan's Liberalized Rule \nD\
\cite{qq:Smullyan:68}.

\PAR{Example}
\label{par:ex:ex3}

The following tree is an \nLKh-proof 
\def\proofSkipAmount{\vskip 0.8ex plus 0.8ex minus 0.4ex}
\begin{prooftree}
\AxiomC{}
\LeftLabel{\nAx}
\UnaryInf$P(a,b) \fCenter ∃y (∃z P(a,z) → P(a,y)), P(a,b)$
\LeftLabel{\nLe}
\RightLabel{$\left(b ≡ \aHc{∃z P(a,z)}\right)$}
\UnaryInf$∃z P(a,z) \fCenter ∃y (∃z P(a,z) → P(a,y)), P(a,b)$
\LeftLabel{\nRi}
\UnaryInf$\fCenter ∃y (∃z P(a,z) → P(a,y)), ∃z P(a,z) → P(a,b)$
\LeftLabel{\nRe}
\UnaryInf$\fCenter ∃y (∃z P(a,z) → P(a,y))$
\LeftLabel{\nRu}
\RightLabel{$\left(a ≡ \aHc{∀x∃y (∃z P(x,z) → P(x,y))}\right)$}
\UnaryInf$\fCenter ∀x ∃y (∃z P(x,z) → P(x,y))$
\end{prooftree}
\def\proofSkipAmount{\vskip 0.8ex plus 0.8ex minus 0.4ex}
of the formula $∀x ∃y (∃z P(x,z) → P(x,y))$.
The Henkin constants $a$ and $b$ belong
to the strong quantifier rule \nRu\ and \nLe, respectively.
Note that the counterexample constant $a$ occurs deeply in the witnessing constant $b$.

When the tree is read as a standard \nLK-proof \nP,
then the proper variable $b$ of the \nLe-inference
violates even Smullyan's liberalized form of eigenvariable condition.
Namely, the principal formula $∃z P(a,z)$ of the \nLe-rule contains
the proper variable of a strong quantifier inference occurring below \nLe;
in this case, it contains the free variable $a$ that belongs to the \nRu-rule.
Even this kind of violation is sound because it satisfies 
the prerequisites of a very liberalized form of eigenvariable condition
in the sequent calculus \nLKpp\ 
(see Par.~\ref{par:lkpp:def} for details).


\section{Properties of \mytexorpdfstring{\nLKh}{LKʰ}-proofs}
\label{sc:pr}
\label{sc:repl}

Throughout the whole section we assume that a first-order language \nLG\ and 
its Henkin expansion \nLGHC\ is fixed.
We also extend the notion of expression to include 
sequents, inference rules, and even proofs.

\PAR{Definition}

We say that two \nLKh-proofs are \emph{structurally similar} if
\begin{itemize}
\item
Their tree representations are isomorphic,
they have the same shape.
\item
Each pair of the corresponding inferences in the isomorphism
can differ only in their side formulas,
both have the same principal and minor formulas.
\end{itemize}
Structurally similar proofs have thus 
the same length, cut-rank and height.

We write $\Gamma ⊆ \Delta$ 
if every element of the multiset $\Gamma$ is 
also an element of the multiset $\Delta$.
Since $A, A, \Gamma ⊆ A, \Gamma$,
the next lemma subsumes not only weakening but contraction as well.

\PAR{Lemma \textnormal{(Weakening, contraction)}}
\label{par:pr:weak}
\label{par:pr:contr}
{\myclaim

Let $\Gamma ⊆ \Pi$ and $\Delta ⊆ \Lambda$ be multisets of formulas.
If \nP\ is an \nLKh-proof of the sequent $\Gamma \fCenter \Delta$,
then we can find a structurally similar \nLKh-proof \nPp\ of the sequent $\Pi \fCenter \Lambda$.

}

\paragraph{Notation.}

We write \aWkP{\nP}{\Pi}{\Lambda} for the proof \nPp\ 
obtained from the proof \nP\ by the application of the lemma.
Our notation is slightly different 
from that of \cite{qq:TroelstraSchwichtenberg:PT:00}.
We thus have
\(
\aCrP{\aWkP{\nP}{\Pi}{\Lambda}} = \aCrP{\nP}
\)
and
\(
\aHtP{\aWkP{\nP}{\Pi}{\Lambda}} = \aHtP{\nP}
\).

\begin{proof}
By a straightforward induction on the structure of the \nLKh-proof \nP.
\end{proof}

\PAR{Definition}

We say that an \nLKh-proof \nP\ is \emph{free for the Henkin constant \aHc{QxA}} 
if there are no strong quantifier rules in \nP\ with $QxA$ as principal formula.

\PAR{Inversion lemma}
\label{par:pr:inv}
{\myclaim

All propositional and strong quantifier \nLKh-rules are invertible 
with length, cut-rank and height preserving inversion.
In particular, we have:

\begin{enumerate}

\item
If
$\nP \rProofLe{l} A → B, \Gamma \fCenter \Delta$,
then we can find
$\nPp_1 \rProofLe{l} \Gamma \fCenter \Delta, A$
and
$\nPp_2 \rProofLe{l} B, \Gamma \fCenter \Delta$
such that
$\aCrP{\nPp_i} ≤ \aCrP{\nP}$ and $\aHtP{\nPp_i} ≤ \aHtP{\nP}$
for $i = 1,2$.

\item
If
$\nP \rProofLe{l} \Gamma \fCenter \Delta, A → B$,
then we can find
$\nPp \rProofLe{l} A, \Gamma \fCenter \Delta, B$
such that
$\aCrP{\nPp} ≤ \aCrP{\nP}$ and $\aHtP{\nPp} ≤ \aHtP{\nP}$.

\item
If
$\nP \rProofLe{l} ∃xA, \Gamma \fCenter \Delta$,
then we can find
$\nPp \rProofLe{l} \aSubst{A}{x}{\aHc{∃xA}}, \Gamma \fCenter \Delta$
free for \aHc{∃xA}
such that
$\aCrP{\nPp} ≤ \aCrP{\nP}$ and $\aHtP{\nPp} ≤ \aHtP{\nP}$.

\end{enumerate}
The similar holds for the remaining inferences.

}

\begin{proof}
This is standard. See, for instance, 
the proof of Lemma 2.5 in \cite{qq:Schwichtenberg:77:proof}.
\end{proof}


\PAR{Deep replacement of terms}
\label{par:repl:def}

We say that a term $r$ \emph{occurs deeply} in an expression $E$ 
if the term $r$  has an (ordinary) occurrence in the expression $E$ or else 
there is a Henkin constant \aHc{∃x A}, 
which has an (ordinary) occurrence in the expression $E$,
such that the term $r$ occurs deeply in its index $∃x A$.
We then say that the term $r$ is a \emph{deep subterm} of the expression $E$.
The subterm is \emph{proper} if $r \not≡ E$.

By a \emph{deep replacement} of a term $r$ by a term $s$ in an expression $E$,
written $E\aReplace{r}{s}$, we mean the replacement of every deep occurrence 
of the term $r$ in the expression $E$ by the term $s$.
The replacement is \emph{non-trivial} if $r \not≡ s$.

\paragraph{Example.}

Consider the Henkin witnessing axiom $∃x P(1,x) → P(1,\aHc{∃x P(1,x)})$.
The constant $1$ has three occurrences in it.
The first two are ordinary occurrences;
the third one --- 
that in the index of the Henkin constant $\aHc{∃x P(1,x)}$ -- 
is deep.
The deep replacement of $1$ by $2$ applied to the axiom 
yields another Henkin witnessing axiom:
\[
\left(∃x P(1,x) → P\!\left(1,\aHc{∃x P(1,x)}\right)\right)\aReplace{1}{2} ≡
\left(∃x P(2,x) → P\!\left(2,\aHc{∃x P(2,x)}\right)\right)
.
\]

\PAR{Inference rules under deep replacement}
\label{par:repl:inf}

We now study the effect of deep replacement on \nLKh-proofs.
This is best explained with an example.
Consider a propositional inference \nLa\ as shown on the left:
\def\proofSkipAmount{\vskip 0.8ex plus 0.8ex minus 0.4ex}
\begin{prooftree}
\AxiomC{$A, B, \GK{A ∧ B,} \Gamma \fCenter \Delta$}
\LeftLabel{\nLa}
\RightLabel{\aReplace{r}{s}}
\UnaryInfC{$A ∧ B, \Gamma \fCenter \Delta$}
\DisplayProof $≡$
\AxiomC{$A\aReplace{r}{s}, B\aReplace{r}{s}, \GK{A\aReplace{r}{s} ∧ B\aReplace{r}{s},} \Gamma\aReplace{r}{s} \fCenter \Delta\aReplace{r}{s}$}
\LeftLabel{\nLa}
\RightLabel{.}
\UnaryInfC{$A\aReplace{r}{s} ∧ B\aReplace{r}{s}, \Gamma\aReplace{r}{s} \fCenter \Delta\aReplace{r}{s}$}
\end{prooftree}
\def\proofSkipAmount{\vskip 0.8ex plus 0.8ex minus 0.4ex}
The effect of the deep replacement \aReplace{r}{s}
on the rule is shown on the right.
The transformation yields another inference of the same kind in the following strong sense:
the principal, minor and side formulas of the first rule
are transformed directly into the principal, minor and side formulas of the second one.
This is because the same changes are done 
both in the premise and in the conclusion of the first inference.

We say that an \nLKh-rule is \emph{invariant} under a deep replacement
if the replacement applied to the rule yields an inference of the same kind
(in the above strong sense).
The effect of deep replacement on invariant rules is shown schematically as follows:
\def\proofSkipAmount{\vskip 0.8ex plus 0.8ex minus 0.4ex}
\begin{prooftree}
\AxiomC{\phantom{$S_1$}}
\LeftLabel{\nRuleAx}
\RightLabel{\aReplace{r}{s}}
\UnaryInfC{$S$} 
\DisplayProof $≡$
\AxiomC{\phantom{$S_1$}}
\LeftLabel{\nRuleAx}
\RightLabel{,}
\UnaryInfC{$S\aReplace{r}{s}$} 
\DisplayProof
\AxiomC{$S_1$}
\LeftLabel{\nRule}
\RightLabel{\aReplace{r}{s}}
\UnaryInfC{$S$} 
\DisplayProof $≡$
\AxiomC{$S_1\aReplace{r}{s}$}
\LeftLabel{\nRule}
\RightLabel{,}
\UnaryInfC{$S\aReplace{r}{s}$} 
\DisplayProof
\quad
\AxiomC{$S_1$}
\AxiomC{$S_2$}
\LeftLabel{\nRule}
\RightLabel{\aReplace{r}{s}}
\BinaryInfC{$S$}
\DisplayProof $≡$
\AxiomC{$S_1\aReplace{r}{s}$}
\AxiomC{$S_2\aReplace{r}{s}$}
\LeftLabel{\nRule}
\RightLabel{.}
\BinaryInfC{$S\aReplace{r}{s}$}
\end{prooftree}
\def\proofSkipAmount{\vskip 0.8ex plus 0.8ex minus 0.4ex}
By the same arguments as for the above \nLa-rule,
it is easy to see that axioms, propositional and cut inferences
are invariant under deep replacements of terms.

This is not so for quantifier inferences.
For an existential quantifier rule 
to be invariant under a deep replacement 
we require that compatible changes are done 
to both its premise and conclusion as shown here:
\def\proofSkipAmount{\vskip 0.8ex plus 0.8ex minus 0.4ex}
\begin{prooftree}
\AxiomC{$\aSubst{A}{x}{\aHc{∃x A}}, \GK{∃x A,} \Gamma \fCenter \Delta$}
\LeftLabel{\nLe}
\RightLabel{\aReplace{r}{s}}
\UnaryInfC{$∃x A, \Gamma \fCenter \Delta$}
\DisplayProof $≡$
\AxiomC{$\aSubstt{A\aReplace{r}{s}}{x}{\aHc{∃x A\aReplace{r}{s}}}, 
         \GK{∃x A\aReplace{r}{s},} \Gamma\aReplace{r}{s} \fCenter 
         \Delta\aReplace{r}{s}$}
\LeftLabel{\nLe}
\RightLabel{,}
\UnaryInfC{$∃x A\aReplace{r}{s}, \Gamma\aReplace{r}{s} \fCenter 
         \Delta\aReplace{r}{s}$}
\end{prooftree}
{\relsize{-0.5}
\def\proofSkipAmount{\vskip 0.8ex plus 0.8ex minus 0.4ex}
\begin{prooftree}
\AxiomC{$\Gamma \fCenter \Delta, ∃x A, \aSubst{A}{x}{t}$}
\LeftLabel{\nRe}
\RightLabel{\aReplace{r}{s}}
\UnaryInfC{$\Gamma \fCenter \Delta, ∃x A$}
\DisplayProof $≡$
\AxiomC{$\Gamma\aReplace{r}{s} \fCenter \Delta\aReplace{r}{s}, ∃x A\aReplace{r}{s},
         \aSubst{A\aReplace{r}{s}}{x}{t\aReplace{r}{s}}$}
\LeftLabel{\nRe}
\RightLabel{.}
\UnaryInfC{$\Gamma\aReplace{r}{s} \fCenter \Delta\aReplace{r}{s}, ∃x A\aReplace{r}{s}$}
\end{prooftree}
\def\proofSkipAmount{\vskip 0.8ex plus 0.8ex minus 0.4ex}
}
Similarly for quantifier inference rules.

\paragraph{Example.}

Consider the following strong quantifier rule shown on the left of the figure
with $∃x P(1,x)$ as principal formula:
\def\proofSkipAmount{\vskip 0.8ex plus 0.8ex minus 0.4ex}
\begin{prooftree}
\AxiomC{$P\!\left(1,\aHc{∃x P(1,x)}\right), \GK{∃x P(1,x),} \Gamma \fCenter \Delta$}
\LeftLabel{\nLe}
\RightLabel{\aReplace{1}{2}}
\UnaryInfC{$∃x P(1,x), \Gamma \fCenter \Delta$}
\DisplayProof $≡$
\AxiomC{$P\!\left(2,\aHc{∃x P(2,x)}\right), \GK{∃x P(2,x),} \Gamma\aReplace{1}{2} \fCenter \Delta\aReplace{1}{2}$}
\LeftLabel{\nLe}
\RightLabel{.}
\UnaryInfC{$∃x P(2,x), \Gamma\aReplace{1}{2} \fCenter \Delta\aReplace{1}{2}$}
\end{prooftree}
\def\proofSkipAmount{\vskip 0.8ex plus 0.8ex minus 0.4ex}
The deep replacement of $1$ by $2$ 
applied to the \nLe-inference 
yields another \nLe-inference with $∃x P(2,x)$ as the principal formula.
The replacement also transforms 
the Henkin constant \aHc{∃x P(1,x)}, which belongs to the first rule,
into the Henkin constant \aHc{∃x P(2,x)}, which belongs to the second one.

This self-correcting nature of deep replacement ---
turning an invariant strong quantifier rule
into another one of the same kind ---
is the reason for extending the customary term replacement 
to operate also on the indices of Henkin constants.

\paragraph{Remark.}

We give here a sufficient condition under which quantifier rules of \nLKh-proofs
are invariant under deep replacement.
Consider a (non-trivial) deep replacement \aReplace{r}{s} in a proof \nP\ 
such that the next two conditions are met:
\begin{enumerate}
\item[(i)]
The term $r$ is different from every Henkin constant 
belonging to some strong quantifier rule applied in \nP.
\item[(ii)]
For every quantifier inference applied in \nP\ 
with the principal formula $Qx A$ and  the minor formula \aSubst{A}{x}{t}, 
there is no subsemiterm $r' \not≡ x$ of $A$,
with at least one free occurrence of the variable $x$,
such that $\aSubst{r'}{x}{t} ≡ r$.
\end{enumerate}
Under 
such conditions every quantifier inference applied in the proof \nP\ 
is invariant under deep replacement of the term $r$ by the term $s$. 

\PAR{Replacement lemma}
\label{par:repl:thm}
{\myclaim

If every rule of the \nLKh-proof
\(
\nP \rProofEq{l} \Gamma \fCenter \Delta
\)
is invariant under deep replacement of a term $r$ by a term $s$, then
\(
\nP\aReplace{r}{s}  \rProofEq{l} 
\Gamma\aReplace{r}{s} \fCenter \Delta\aReplace{r}{s}
\).

}

\begin{proof}
By a straightforward induction on the structure of the \nLKh-proof \nP.
\end{proof}


\section{Cut Elimination}
\label{sc:cut}

The cut elimination theorem for the sequent calculus \nLKh\ states that
\begin{quote}
{\myclaim

if \nP\ is an \nLKh-proof of a pure formula,
then there is a cut-free \nLKh-proof \nPp\ of the same formula 
with height
\(
\aHtP{\nPp} ≤ \aSExp{2\aCrP{\nP}}{\aHtP{\nP}}
\).

}
\end{quote}
The proof of the theorem cannot rely on the Gentzen-style of cut elimination,
because \nLKh-proofs are not \emph{regular},
they do not satisfy the standard eigenvariable condition.
Regularization of an \nLKh-proof by the method 
described in Komara and Voda \cite[Thm.~15]{qq:KomaraVoda:95:TABLEAUX95}
does not help either since such regularization introduces new cuts 
that cannot be expressed in term
of the cut-rank of the original (non-regular) proof.
The proposed cut elimination for \nLKh\ intermingles two steps: 
quasiregularization and the standard cut-rank reduction.
Quasiregularization restores a weaker form of regularity
for quantifier cuts with the highest rank,
it preserves proof cut-rank but increases proof height exponentially.
Compare with Aguilera and Baaz~\cite[Cor. 4.3]{DBLP:journals/jsyml/AguileraB19}.

\PAR{Quasiregularity condition}

We say that a cut rule of the form ($Q ∈ \{∃,∀\}$)
\def\proofSkipAmount{\vskip 0.8ex plus 0.8ex minus 0.4ex}
\begin{prooftree}
\AxiomC{$\Gamma \fCenter \Delta, QxA$}
\AxiomC{$QxA, \Gamma \fCenter \Delta$}
\LeftLabel{\nCut}
\RightLabel{}
\BinaryInfC{$\Gamma \fCenter \Delta$}
\end{prooftree}
\def\proofSkipAmount{\vskip 0.8ex plus 0.8ex minus 0.4ex}
is \emph{critical} if the Henkin constant \aHc{QxA},
which belongs to the quantifier formula $QxA$,
deeply occurs in the conclusion $\Gamma \fCenter \Delta$ of the cut inference.

An \nLKh-proof is \emph{$r$-quasiregular}
if it does not contain any critical cut of rank $r$.
Quasiregularity is downwards hereditary property:
subderivations of an $r$-quasi\-reg\-u\-lar proof are $r$-quasiregular as well.
This holds also in reverse even for the case 
when a proof ends with a cut of rank $r$ on a quantifier formula provided that 
the Henkin constant that belongs to the cut-formula does not deeply occur in the endsequent.

\PAR{Inversion lemma for cut rules}
\label{par:cut:inv:cut}
{\myclaim

Let
\(
\nP \rProof \Gamma \fCenter \Delta
\)
be an \nLKh-proof and $A$ a formula.
Then we can find \nLKh-proofs
\(
\nP[\rho] \rProof \Gamma \fCenter \Delta, A
\)
and
\(
\nP[\sigma] \rProof A, \Gamma \fCenter \Delta
\)
with height $\aHtP{\nP_i} ≤ \aHtP{\nP}$ for $i = 1,2$
that contain together, apart from $A$, the same cut-formulas as \nP.

}

\begin{proof}
This is done by pushing every cut on $A$ in \aP{} downwards over the remaining inferences.
Merging consecutive cuts on $A$ into a single one requires contraction.
\end{proof}

\PAR{Quasiregularization lemma}
\label{par:cut:qreg}
{\myclaim

Let \nP\ be an \nLKh-proof of a pure sequent. 
Then we can find an \aCrP{\nP}-quasiregular \nLKh-proof \nPp\ with the same endsequent 
such that $\aCrP{\nPp} = \aCrP{\nP}$ and $\aHtP{\nPp} < 2^{\aHtP{\nP}}$.

}

\begin{proof}

By a repeated application of the cut-inversion lemma
to quantifier cut-for\-mu\-las of depth \aCrP{\nP} in the proof \nP.
For that, we first prove the following claim:
\begin{quote}
{\myclaim

Let $n$ be the number 
of different quantifier cut-formulas  of depth $r$ occurring in the \nLKh-proof
\(
\nP \rProof \Gamma \fCenter \Delta
\)
with cut-rank $\aCrP{\nP} ≤ r$.
Let further there is no Henkin constant deeply 
occurring in the sequent $\Gamma \fCenter \Delta$ 
that belongs to a quantifier cut-formula of depth $r$ in \nP.
Then there is an $r$-quasiregular \nLKh-proof 
\(
\nPp \rProof \Gamma \fCenter \Delta
\)
containing the same quantifier cut-formulas of depth $r$ as \nP\
with cut-rank $\aCrP{\nPp} = \aCrP{\nP}$ and 
with height $\aHtP{\nPp} ≤ \aHtP{\nP}+n$.

}
\end{quote}
This is proved by induction on the number $n$.
If $n=0$, then the proof \nP\ is already $r$-quasiregular and 
it suffices to take \nP\ for \nPp.
If $n > 0$, then the derivation \nP\ contains at least one quantifier cut-formula $A$ of depth $r$.
We may suppose without loss of generality that the Henkin constant \aHc{A}
does not deeply occur in another quantifier cut-formula of the same depth $r$.
By the cut inversion lemma~\ref{par:cut:inv:cut} we can find proofs
\begin{align*}
\nP_1 \rProof \Gamma \fCenter \Delta, A
\quad \text{and} \quad
\nP_2 \rProof A, \Gamma \fCenter \Delta
\end{align*}
that contain together, apart from $A$, the same quantifier cut-formulas of depth $r$ 
as the proof \nP.
Also $\aCrP{\nP_i} ≤ \aCrP{\nP}$ and $\aHtP{\nP_i} ≤ \aHtP{\nP}$ for $i = 1,2$ 
by the same lemma.
We apply the inductive hypothesis to both derivations \aP{1}, \aP{2}
and obtain $r$-quasiregular proofs
\begin{align*}
\aPp{1} \rProof \Gamma \fCenter \Delta, A
\quad \text{and} \quad
\aPp{2} \rProof A, \Gamma \fCenter \Delta
\end{align*}
that contain together the same quantifier cut-formulas of depth $r$ 
as the derivations \aP{1} and \aP{2}.
Also $\aCrP{\aPp{i}} = \aCrP{\aP{i}}$ and $\aHtP{\aPp{i}} ≤ \aHtP{\aP{i}}+n-1$ 
by IH for $i = 1,2$.
We built a new deduction ending with a cut inference on $A$:
\def\proofSkipAmount{\vskip 0.8ex plus 0.8ex minus 0.4ex}
\begin{prooftree}
\AxiomC{}
\RightLabel{\aPp{1}}
\branchDeduce
\DeduceC{$\Gamma \fCenter \Delta, A$}
\AxiomC{}
\RightLabel{\aPp{2}}
\branchDeduce
\DeduceC{$A, \Gamma \fCenter \Delta$}
\LeftLabel{\nCut}
\RightLabel{.}
\BinaryInfC{$\Gamma \fCenter \Delta$}
\RightSubproofLabel{\nPp}
\end{prooftree}
\def\proofSkipAmount{\vskip 0.8ex plus 0.8ex minus 0.4ex}
By assumption, the Henkin constant \aHc{A} does not occur deeply in $\Gamma \fCenter \Delta$.
By IH, both derivations \aPp{1} and \aPp{2} are $r$-quasiregular.
Hence so is the proof \nPp. We also have
\begin{align*}
\aCrP{\nPp} & =
\max(\aDpF{A},\aCrP{\aPp{1}},\aCrP{\aPp{2}}) \mEQx{$\aDpF{A} = \aCrP{\nP}$}{=} 
\max(\aCrP{\nP},\aCrP{\aPp{1}},\aCrP{\aPp{2}}) 
\\ & \mEQx{IH}{=} 
\max(\aCrP{\nP},\aCrP{\aP{1}},\aCrP{\aP{2}}) \mEQx{L \ref{par:cut:inv:cut}}{=} 
\aCrP{\nP},
\\
\aHtP{\nPp} & =
\max\left(\aHtP{\aPp{1}},\aHtP{\aPp{2}}\right)+1 =
\max\left(\aHtP{\aPp{1}}+1,\aHtP{\aPp{2}}+1\right) 
\\ & \mEQx{IH}{≤}
\max\left(\aHtP{\aP{1}}+n,\aHtP{\aP{2}}+n\right) =
\max\left(\aHtP{\aP{1}},\aHtP{\aP{2}}\right)+n \mEQx{L \ref{par:cut:inv:cut}}{≤}
\aHtP{\nP}+n
.
\end{align*}
This ends the proof of the auxiliary claim.

We can now prove the lemma.
Let \nP\ be a proof with a pure endsequent and
let $n$ be the number of different quantifier cut-formulas 
of depth \aCrP{\nP} occurring in the proof \nP.
By the auxiliary claim
we can find an \aCrP{\nP}-quasiregular derivation \nPp\ with the same endsequent 
such that $\aCrP{\nPp} = \aCrP{\nP}$ and $\aHtP{\nPp} ≤ \aHtP{\nP}+n$.
Because the proof \nP\ contains at least $n$ cuts, 
it must be $2n < 2^{\aHtP{\nP}}$.
Hence
\(
n < 2^{\aHtP{\nP}-1} ≤ 2^{\aHtP{\nP}} -\aHtP{\nP}
\)
and thus
\(
\aHtP{\nPp} ≤ \aHtP{\nP}+n < \aHtP{\nP}+2^{\aHtP{\nP}} -\aHtP{\nP} = 2^{\aHtP{\nP}}
\).
\end{proof}

\PAR{Lemma}
\label{par:cut:red:I:e}
{\myclaim

Let $\Gamma ⊆ \Pi$ and $\Delta ⊆ \Lambda$ be multisets of formulas.
Let further the Henkin constant \aHc{∃xA} does not occur deeply 
in the sequent $\Gamma \fCenter \Delta$ and let $r = \aDpF{∃xA}$.
Let finally
\(
\aP{1} \rProof \Pi \fCenter \Lambda, ∃xA
\)
and
\(
\aP{2} \rProof \aSubst{A}{x}{\aHc{∃xA}}, \Gamma \fCenter \Delta
\)
be \nLKh-proofs with cuts of rank $< r$ such that the proof \aP{2} is free for \aHc{∃xA}.
Then there is an \nLKh-proof 
\(
\nPp \rProof \Pi \fCenter \Lambda
\)
with all its cuts of rank $< r$ and with height 
\(
\aHtP{\nPp} ≤ \aHtP{\aP{1}}+\aHtP{\aP{2}}
\).

}

\begin{proof}

By induction on the structure of the \nLKh-proof \aP{1}.
We continue by case analysis of the last rule applied in \aP{1}.

Suppose that the proof \aP{1} ends with an \nRe-inference of the form
\def\proofSkipAmount{\vskip 0.8ex plus 0.8ex minus 0.4ex}
\begin{prooftree}
\AxiomC{}
\RightLabel{\aP{11}}
\branchDeduce
\DeduceC{$\Pi \fCenter \Lambda, ∃xA, \aSubst{A}{x}{t}$}
\LeftLabel{\nRe}
\RightLabel{.}
\UnaryInfC{$\Pi \fCenter \Lambda, ∃xA$}
\RightSubproofLabel{\aP{1}}
\end{prooftree}
\def\proofSkipAmount{\vskip 0.8ex plus 0.8ex minus 0.4ex}
We apply the inductive hypothesis to \aP{11} and find a derivation
\(
\aPp{1} \rProof \Pi \fCenter \Lambda, \aSubst{A}{x}{t}
\)
with cuts of rank $< r$ and height 
\(
\aHtP{\aPp{1}} ≤ \aHtP{\aP{11}}+\aHtP{\aP{2}}
\).
The proof \aP{2} is free for \aHc{∃xA} and hence
\begin{gather*}
\aWkP{\aP{2}\aReplace{\aHc{∃xA}}{t}}{\aSubst{A}{x}{t}, \Pi}{\Lambda}
\rProof \aSubst{A}{x}{t}, \Pi \fCenter \Lambda
,
\\
\aCrP{\aWkP{\aP{2}\aReplace{\aHc{∃xA}}{t}}{\aSubst{A}{x}{t}, \Pi}{\Lambda}} =
\aCrP{\aP{2}} < r
,
\\
\aHtP{\aWkP{\aP{2}\aReplace{\aHc{∃xA}}{t}}{\aSubst{A}{x}{t}, \Pi}{\Lambda}} =
\aHtP{\aP{2}}
\end{gather*}
by the replacement lemma \ref{par:repl:thm} and the weakening lemma \ref{par:pr:weak}.
We now use a cut on the formula $\aSubst{A}{x}{t}$ and obtain a new derivation
\def\proofSkipAmount{\vskip 0.8ex plus 0.8ex minus 0.4ex}
\begin{prooftree}
\AxiomC{}
\RightLabel{\aPp{1}}
\branchDeduce
\DeduceC{$\Pi \fCenter \Lambda, \aSubst{A}{x}{t}$}
\AxiomC{}
\RightLabel{\aWkP{\aP{2}\aReplace{\aHc{∃xA}}{t}}{\aSubst{A}{x}{t}, \Pi}{\Lambda}}
\branchDeduce
\DeduceC{$\aSubst{A}{x}{t}, \Pi \fCenter \Lambda$}
\LeftLabel{\nCut}
\RightLabel{.}
\BinaryInfC{$\Pi \fCenter \Lambda$}
\RightSubproofLabel{\nPp}
\end{prooftree}
\def\proofSkipAmount{\vskip 0.8ex plus 0.8ex minus 0.4ex}
The proof \nPp\ has all its cuts with rank $< r$
since $\aDpF{\aSubst{A}{x}{t}} < r$.
Its height satisfies 
\begin{align*}
\aHtP{\nPp} & = 
\max\left(
\aHtP{\aPp{1}},
\aHtP{\aWkP{\aP{2}\aReplace{\aHc{∃xA}}{t}}{\aSubst{A}{x}{t}, \Pi}{\Lambda}}
\right)+1 =
\max\left(\aHtP{\aPp{1}},\aHtP{\aP{2}}\right)+1 
\\ & \mEQx{IH}{≤}
\max\left(\aHtP{\aP{11}}+\aHtP{\aP{2}},\aHtP{\aP{2}}\right)+1 =
\aHtP{\aP{11}}+\aHtP{\aP{2}}+1 =
\aHtP{\aP{1}}+\aHtP{\aP{2}}
.
\end{align*}

The case when the proof \aP{1} ends with an axiom is obvious
because existential formulas are non-atomic.
The proof of the remaining cases is a straightforward application of the inductive hypothesis.
\end{proof}

\PAR{Reduction lemma I}
\label{par:cut:red:I}
{\myclaim

Let 
\(
\nP \rProof \Gamma \fCenter \Delta
\)
be a $r$-quasiregular \nLKh-proof with a final inference a cut of rank $r$ 
such that every other cut in \nP\ has rank $< r$.
Then there is an \nLKh-proof 
\(
\nPp \rProof \Gamma \fCenter \Delta
\)
with all its cuts of rank $< r$ and with height $\aHtP{\nPp} < 2\aHtP{\nP}$.

}

\begin{proof}

The $r$-quasiregular proof \nP\ is of the form
\def\proofSkipAmount{\vskip 0.8ex plus 0.8ex minus 0.4ex}
\begin{prooftree}
\AxiomC{}
\RightLabel{\aP{1}}
\branchDeduce
\DeduceC{$\Gamma \fCenter \Delta, A$}
\AxiomC{}
\RightLabel{\aP{2}}
\branchDeduce
\DeduceC{$A, \Gamma \fCenter \Delta$}
\LeftLabel{\nCut}
\RightLabel{}
\BinaryInfC{$\Gamma \fCenter \Delta$}
\RightSubproofLabel{\nP}
\end{prooftree}
\def\proofSkipAmount{\vskip 0.8ex plus 0.8ex minus 0.4ex}
with cut-rank $r = \aDpF{A}$ and $\aCrP{\aP{i}} < r$ for $i=1,2$.
The proof proceeds by case analysis according to the form of the cut-formula $A$.

\emph{Case:}
$A ≡ ∃xB$.
By the inversion lemma \ref{par:pr:inv} we construct a new proof
\def\proofSkipAmount{\vskip 0.8ex plus 0.8ex minus 0.4ex}
\begin{prooftree}
\AxiomC{}
\RightLabel{\aP{1}}
\branchDeduce
\DeduceC{$\Gamma \fCenter \Delta, ∃xB$}
\AxiomC{}
\RightLabel{\aPp{2}}
\branchDeduce
\DeduceC{$\aSubst{B}{x}{\aHc{∃xB}}, \Gamma \fCenter \Delta$}
\LeftLabel{\nLe}
\RightLabel{}
\UnaryInfC{$∃xB, \Gamma \fCenter \Delta$}
\LeftLabel{\nCut}
\RightLabel{,}
\BinaryInfC{$\Gamma \fCenter \Delta$}
\RightSubproofLabel{\nPpp}
\end{prooftree}
\def\proofSkipAmount{\vskip 0.8ex plus 0.8ex minus 0.4ex}
where its subderivation \aPp{2} is free for \aHc{∃xB} 
with cut-rank
\(
\aCrP{\aPp{2}} < r
\)
and with height
\(
\aHtP{\aPp{2}} ≤ \aHtP{\aP{2}}
\).
Moreover, from the $r$-quasiregularity of \nP\ we know that
the Henkin witnessing constant \aHc{∃xB} does not occur deeply 
in the sequent $\Gamma \fCenter \Delta$.
\emph{This is the place in the proof, where the quasiregularity assumption is crucial.}
By Lemma~\ref{par:cut:red:I:e} 
applied to \aP{1} and \aPp{2}
we can find a proof 
\(
\nPp \rProof \Gamma \fCenter \Delta
\)
with cut-rank $< r$ and with height
\begin{align*}
\aHtP{\nPp} \mEQx{L \ref{par:cut:red:I:e}}{≤} 
\aHtP{\aP{1}} + \aHtP{\aPp{2}}  \mEQx{L \ref{par:pr:inv}}{≤}
\aHtP{\aP{1}} + \aHtP{\aP{2}} <
2 (\max(\aHtP{\aP{1}},\aHtP{\aP{2}}) + 1) =
2\aHtP{\nP}
.
\end{align*}
This proves the case when $A$ is an existential formula.

The case when $A$ is a universal formula is dual to the previous one.
For the proof of the remaining cases see, for instance, 
the proof of Lemma 2.4.2.1 in \cite{qq:Buss:ipt:98} or, alternatively, 
the proof of Refined Reduction Lemma in \cite{qq:Gerhardy:03:CSL03}.
\end{proof}

\PAR{Reduction lemma II}
\label{par:cut:red:II}
{\myclaim

Let
\(
\nP \rProof \Gamma \fCenter \Delta
\)
be a $r$-quasiregular \nLKh-proof with cuts of rank $≤ r$.
Then there is an \nLKh-proof 
\(
\nPp \rProof \Gamma \fCenter \Delta
\)
with cuts of rank $< r$ and with height $\aHtP{\nPp} < 2^{\aHtP{\nP}}$.

}

\begin{proof}
By induction on the structure of the \nLKh-proof \nP.
We distinguish several cases according to the last rule applied in \nP.

Suppose that the $r$-quasiregular proof \nP\ ends with a cut inference of the form
\def\proofSkipAmount{\vskip 0.8ex plus 0.8ex minus 0.4ex}
\begin{prooftree}
\AxiomC{}
\RightLabel{\aP{1}}
\branchDeduce
\DeduceC{$\Gamma \fCenter \Delta, A$}
\AxiomC{}
\RightLabel{\aP{2}}
\branchDeduce
\DeduceC{$A, \Gamma \fCenter \Delta$}
\LeftLabel{\nCut}
\RightLabel{.}
\BinaryInfC{$\Gamma \fCenter \Delta$}
\RightSubproofLabel{\nP}
\end{prooftree}
\def\proofSkipAmount{\vskip 0.8ex plus 0.8ex minus 0.4ex}
Here
\(
\aP{1} \rProof \Gamma \fCenter \Delta, A
\)
and
\(
\aP{2} \rProof A, \Gamma \fCenter \Delta
\)
are $r$-quasiregular proofs with all their cuts of rank $≤ r$.
We apply the inductive hypothesis to both subderivations and obtain proofs
\(
\aPp{1} \rProof \Gamma \fCenter \Delta, A
\)
and
\(
\aPp{2} \rProof A, \Gamma \fCenter \Delta
\)
with all their cuts of rank $< r$ and 
with height $\aHtP{\aPp{i}} < 2^{\aHtP{\aP{i}}}$ for $i=1,2$.
We use a similar cut inference and built a new proof
\def\proofSkipAmount{\vskip 0.8ex plus 0.8ex minus 0.4ex}
\begin{prooftree}
\AxiomC{}
\RightLabel{\aPp{1}}
\branchDeduce
\DeduceC{$\Gamma \fCenter \Delta, A$}
\AxiomC{}
\RightLabel{\aPp{2}}
\branchDeduce
\DeduceC{$A, \Gamma \fCenter \Delta$}
\LeftLabel{\nCut}
\RightLabel{.}
\BinaryInfC{$\Gamma \fCenter \Delta$}
\RightSubproofLabel{\nPpp}
\end{prooftree}
\def\proofSkipAmount{\vskip 0.8ex plus 0.8ex minus 0.4ex}
Note that we have
\begin{align}
\aHtP{\nPpp} =
\max\left(\aHtP{\aPp{1}},\aHtP{\aPp{2}}\right)+1 \mEQx{IH}{≤}
\max\left(2^{\aHtP{\aP{1}}},2^{\aHtP{\aP{2}}}\right) =
2^{\max(\aHtP{\aP{1}},\aHtP{\aP{2}})}
.
\label{eq:cut:cut1:qreg:Ppp}
\end{align}
We now consider two subcases according to the depth \aDpF{A} of the cut-formula $A$.

Subcase: $\aDpF{A} = r$.
The proof \nPpp\ is $r$-quasiregular,
for otherwise $A$ would be a quantifier formula such that
its Henkin constant \aHc{A} would deeply occur in $\Gamma \fCenter \Delta$,
but the last is impossible due to $r$-quasiregularity of \nP.
By the reduction lemma~\ref{par:cut:red:I} applied to \nPpp,
we can find a proof 
\(
\nPp \rProof \Gamma \fCenter \Delta
\)
with cut-rank $< r$ and with height
\begin{align*}
\aHtP{\nPp} & 
\mEQx{L \ref{par:cut:red:I}}{<}
2\aHtP{\nPpp} \mEQx{\XEL{cut:cut1:qreg}{Ppp}}{≤}
2 \times 2^{\max(\aHtP{\aP{1}},\aHtP{\aP{2}})} =
2^{\max(\aHtP{\aP{1}},\aHtP{\aP{2}})+1} =
2^{\aHtP{\nP}}
.
\end{align*}
Subcase: $\aDpF{A} < r$.
Then every cut in the proof \nPpp\ has rank $< r$ and 
\begin{align*}
\aHtP{\nPpp} \mEQx{\XEL{cut:cut1:qreg}{Ppp}}{≤}
2^{\max(\aHtP{\aP{1}},\aHtP{\aP{2}})} <
2^{\max(\aHtP{\aP{1}},\aHtP{\aP{2}})+1} =
2^{\aHtP{\nP}}
.
\end{align*}
So it suffices to take \nPpp\ for \nPp.

The case when the proof \aP{} ends with an axiom is easy to check. 
The proof of the remaining cases is a straightforward application of the inductive hypothesis.
\end{proof}

\PAR{Reduction lemma III}
\label{par:cut:red:III}
{\myclaim

Let \nP\ be an \nLKh-proof of a pure sequent with cut-rank $\aCrP{\nP} > 0$.
Then there exists an \nLKh-proof \nPp\ with the same endsequent such that
$\aCrP{\nPp} < \aCrP{\nP}$ and $\aHtP{\nPp} < 2^{2^{\aHtP{\nP}}}$.

}

\begin{proof}
By the quasiregularization lemma~\ref{par:cut:qreg} and 
the reduction lemma~\ref{par:cut:red:II}.
\end{proof}

\PAR{Cut elimination theorem}
\label{par:cut:thm}
{\myclaim

Let \nP\ be an \nLKh-proof of a pure sequent.
Then there is a cut-free \nLKh-proof \nPp\ with the same endsequent such that
\(
\aHtP{\nPp} ≤ \aSExp{2\aCrP{\nP}}{\aHtP{\nP}}
\).

}

\begin{proof}
By repeated application of the reduction lemma~\ref{par:cut:red:III}.
\end{proof}

\PAR{Theorem \textnormal{(Subformula property)}}
\label{par:cut:sub}
{\myclaim

If \nP\ is a cut-free \nLKh-proof, 
then every formula of \nP\ is a subformula 
of some formula from the endsequent of \nP.

}

\begin{proof}
By a straightforward induction on the structure of the \nLKh-proof \nP.
\end{proof}

\PAR{Example \textnormal{(Komara and Voda \mytexorpdfstring{\cite{qq:KomaraVoda:95:TABLEAUX95}}{1995}})}
\label{par:cut:ex}

Cut elimination does not extend to \nLKh-proofs with non-pure endsequents.
Indeed, consider the \nLKh-proof with a cut
\def\proofSkipAmount{\vskip 0.8ex plus 0.8ex minus 0.4ex}
\begin{prooftree}\GK{\small}
\AxiomC{}
\LeftLabel{\nAx}
\UnaryInf$P(k) \fCenter \GK{P(k) → P(a), P(a),} ∃x P(x), P(k)$
\LeftLabel{\nRe}
\UnaryInf$P(k) \fCenter \GK{P(k) → P(a), P(a),} ∃x P(x)$
\AxiomC{}
\LeftLabel{\nAx}
\UnaryInf$P(a)\GK{, ∃x P(x), P(k)} \fCenter \GK{P(k) → P(a),} P(a)$
\LeftLabel{\nLe}
\RightLabel{$\left(a ≡ \aHc{∃xP(x)}\right)$}
\UnaryInf$∃x P(x)\GK{, P(k)} \fCenter \GK{P(k) → P(a),} P(a)$
\LeftLabel{\nCut}
\BinaryInf$P(k) \fCenter \GK{P(k) → P(a),} P(a)$
\LeftLabel{\nRi}
\RightLabel{}
\UnaryInf$\fCenter P(k) → P(a)$
\end{prooftree}
\def\proofSkipAmount{\vskip 0.8ex plus 0.8ex minus 0.4ex}
of the (non-pure) formula
\begin{align}
P(k) → P(a).
\label{eq:cut:ex:f}
\end{align}
The $k$ is a pure constant and $a ≡ \aHc{∃xP(x)}$ 
is the Henkin witnessing constant for $∃xP(x)$.
We claim that the formula \XEL{cut:ex}{f} cannot have a cut-free \nLKh-proof.
Assume, by contradiction, that it has such a proof.
By the subformula property, the proof consists entirely 
of the subformulas of the formula \XEL{cut:ex}{f}.
Therefore it does not contain any quantifier inferences. 
In fact, it is a propositional proof.
But this is impossible, 
because the formula \XEL{cut:ex}{f} is certainly not a tautology.


\section{Eliminating Skolem functions}
\label{sc:skk}

Maehara's method \cite[Lemma 8.11]{qq:Takeuti:87}
for eliminating Skolem functions from cut-free 
\nLK-derivations is based on the following simple idea.
Let \nP\ be a cut-free proof of a formula $∀x A[x,f(x)] → B$ with length $l$,
where $f$ is a new function symbol.
Find all \emph{Skolem terms}
\(
f(s_1), \ldots, f(s_n)
\)
in \nP\ that belong to \nLu-rules with $∀x A[x,f(x)]$ as principal formula.
Replace them with fresh variables \aSeq{n}{a} and then
quantify them with \nLe-rules taking
\(
∃y A[s_1,y], \ldots, ∃y A[s_n,y]
\)
as principal formulas.
The trouble is that the resulting proof-like structure \nPp\ is not a correct \nLK-proof,
for the variables \aSeq{n}{a} may violated the eigenvariable condition.
So the tree \nPp\ is further transformed by introducing cuts to restore the eigenvariable condition.
This yields a correct \nLK-proof of the formula $∀x ∃y A[x,y] → B$.
The subsequent cut elimination leads to superexponential blow-up in proof length.

For \nLKh-proofs the step needed to restore eigenvariable condition is not necessary.
The \nPp\ is already a valid \nLKh-proof with length $l+n$ if the transformation is done properly.
To see this, consider a cut-free \nLKh-proof \nP\ with length $l$ and $n=2$:
\def\proofSkipAmount{\vskip -1.0ex plus 0.8ex minus 0.4ex}
{\relsize{-0.5}
\begin{prooftree}
\hskip -4.0em
\AxiomC{}
\RightLabel{\aP[\rho]{}}
\branchDeduce
\DeduceC{$A[r,f(r)], ∀x A[x,f(x)], \Gamma \fCenter \Delta$}
\LeftLabel{\nLu}
\UnaryInfC{$∀x A[x,f(x)], \Gamma \fCenter \Delta$}
\AxiomC{}
\RightLabel{\aP[\sigma]{}}
\branchDeduce
\DeduceC{$A[s,f(s)], ∀x A[x,f(x)], \Pi \fCenter \Lambda$}
\LeftLabel{\nLu}
\UnaryInfC{$∀x A[x,f(x)], \Pi \fCenter \Lambda$}
\noLine
\def\extraVskip{-1pt}
\BinaryInfC{}
\def\extraVskip{2pt}
\branchDeduce
\DeduceC{$∀x A[x,f(x)] \fCenter B$.}
\RightSubproofLabel{\nP}
\end{prooftree}
\def\proofSkipAmount{\vskip 0.8ex plus 0.8ex minus 0.4ex}
}
\noindent
These two \nLu-rules are the only inferences in \aP{} that are
non-invariant under deep replacement of $f$-terms.
We may suppose without loss of generality that 
the Skolem term $f(r)$ is not a deep subterm of the Skolem term $f(s)$. 
We eliminate $f(r)$ and $f(s)$ by deeply replacing them, as described below, 
with the Henkin constants $\aHc{∃y A[r,y]}$ and $\aHc{∃y A[s,y]}$.
Let us denote by \nDr\ the following composition of deep replacements:
\[
\nDr ≡ \aReplacee{f(r)}{\aHc{∃y A[r,y]}}\aReplacee{f(s)}{\aHc{∃y A[s,y]}}
.
\]
From
\(
f(s)\aReplace{f(r)}{\aHc{∃y A[r,y]}} ≡ f(s)
\)
we obtain that
\begin{align*}
\aSubst{A}{}{r,f(r)}\nDr ≡
\aSubstt{A}{}{r\nDr,\aHc{∃y A[r\nDr,y]}} 
,
\qquad
\aSubst{A}{}{s,f(s)}\nDr ≡
\aSubstt{A}{}{s\nDr,\aHc{∃y A[s\nDr,y]}}
.
\end{align*}
By the replacement lemma~\ref{par:repl:thm},
the following is a valid \nLKh-proof of length $l+2$:
\def\proofSkipAmount{\vskip -1.0ex plus 0.8ex minus 0.4ex}
{\relsize{-1.0}
\begin{prooftree}
\hskip -5.2em
\AxiomC{}
\RightLabel{\aP[\rho]{}\nDr}
\branchDeduce
\DeduceC{$\aSubstt{A}{}{r\nDr,\aHc{∃y A[r\nDr,y]}}, ∀x∃y A[x,y], \Gamma\nDr \fCenter \Delta\nDr$}
\LeftLabel{\nLe}
\UnaryInfC{$∃y A[r\nDr,y], ∀x∃y A[x,y], \Gamma\nDr \fCenter \Delta\nDr$}
\LeftLabel{\nLu}
\UnaryInfC{$∀x∃y A[x,y], \Gamma\nDr \fCenter \Delta\nDr$}
\AxiomC{}
\RightLabel{\aP[\sigma]{}\nDr}
\branchDeduce
\DeduceC{$\aSubstt{A}{}{s\nDr,\aHc{∃y A[s\nDr,y]}}, ∀x∃y A[x,y], \Pi\nDr \fCenter \Lambda\nDr$}
\LeftLabel{\nLe}
\UnaryInfC{$∃yA[s\nDr,y], ∀x∃y A[x,y], \Pi\nDr \fCenter \Lambda\nDr$}
\LeftLabel{\nLu}
\UnaryInfC{$∀x∃y A[x,y], \Pi\nDr \fCenter \Lambda\nDr$}
\noLine
\def\extraVskip{-1pt}
\BinaryInfC{}
\def\extraVskip{2pt}
\branchDeduce
\RightLabel{$\nDr$}
\DeduceC{$∀x∃y A[x,y] \fCenter B$.}
\RightSubproofLabel{\nPp}
\end{prooftree}
\def\proofSkipAmount{\vskip 0.8ex plus 0.8ex minus 0.4ex}
}
\noindent
A similar idea in Aguilera and Baaz \cite[Prop. 4.8]{DBLP:journals/jsyml/AguileraB19},
but without a detailed proof, is used to show that deskolemization
of cut-free \nLKpp-proofs is linear.

\PAR{Skolem witnessing functions}
\label{par:sk:wf:def}

Let \nLG\ be a first-order language and $f$ a new $k$-ary function symbol not \nLG.
Let further
\begin{align}
\aSeqc{k}{∀x}∃yA[\aSeq{k}{x},y]
\label{eq:sk:wf:def:f}
\end{align}
be an \nLG-formula. 
A \emph{one-step witnessing Skolemization} of 
\XEL{sk:wf:def}{f} 
is a formula
\begin{align}
\aSeqc{k}{∀x} A[\aSeq{k}{x}, f(\aSeq{k}{x})]
\label{eq:sk:wf:def:sk}
.
\end{align}
The function $f$ is called the \emph{Skolem witnessing function} of the skolemization.

We say that an \nLKh-proof \nP\ is 
\emph{free for the Skolem witnessing function $f$} 
(\emph{free for $f$} for short)
if the only rules in the proof \nP\
that are non-invariant under deep replacement of $f$-terms
are \nLu-inferences of the form
{\relsize{-0.5}
\def\proofSkipAmount{\vskip 0.8ex plus 0.8ex minus 0.4ex}
\begin{prooftree}
\AxiomC{$
A[\aSeqnq{k-1}{s}{i}, s_k, f(\aSeqnq{k-1}{s}{i}, s_k)], 
∀x_k A[\aSeqnq{k-1}{s}{i}, x_k, f(\aSeqnq{k-1}{s}{i}, x_k)], 
\Pi \fCenter \Lambda
$}
\LeftLabel{\nLu}
\RightLabel{.}
\UnaryInfC{$
∀x_k A[\aSeqnq{k-1}{s}{i}, x_k, f(\aSeqnq{k-1}{s}{i}, x_k)], 
\Pi \fCenter \Lambda
$}
\end{prooftree}
\def\proofSkipAmount{\vskip 0.8ex plus 0.8ex minus 0.4ex}
}
\noindent
By the subformula property~\ref{par:cut:sub} cut-free proofs are trivially free for $f$.
The property of an \nLKh-proof being free for $f$ is downwards hereditary property:
subderivations of such a proof are themselves free for $f$.

In the next theorem and its corollary we assume that 
the first-order language \nLG, 
the Skolem witnessing function $f$ and
the skolemization formula 
are the same as above.
Note that we are working in the Henkin expansion 
of the first-order language $\nLG \cup \{f\}$.

\PAR{Elimination theorem \textnormal{(Skolem witnessing function)}}
\label{par:sk:wf:thm}
{\myclaim

Let $\Gamma$ and $\Delta$ be finite multisets of \nLG-formulas,
Let 
further
\[
\nP \rProofEq{l} \aSeqc{k}{∀x} A[\aSeq{k}{x},f(\aSeq{k}{x})], \Gamma \fCenter \Delta
\]
be an \nLKh-proof free for the Skolem witnessing function $f$. 
Let finally $n$ be the number of inferences in \aP{} 
which are non-invariant under deep replacement of $f$-terms.
Then there is an \nLKh-proof \nPp\ in \nLG\ 
with the same cut-rank as \nP\ such that
\[
\nPp \rProofEq{l+n} \aSeqc{k}{∀x} ∃y A[\aSeq{k}{x},y], \Gamma \fCenter \Delta
.
\]

}

\begin{proof}
\renewcommand{\qedsymbol}{\ensuremath{\blacksquare}}

By assumption, the only rules in the proof \aP{} 
that are non-invariant under deep replacement of $f$-terms
are \nLu-inferences of the form ($i = 1, \ldots, n$)
{\relsize{-0.5}
\def\proofSkipAmount{\vskip -2.0ex plus 0.8ex minus 0.4ex}
\begin{prooftree}
\AxiomC{$
A[\aSeqq{k-1}{s}{i}, s_k^i, f(\aSeqq{k-1}{s}{i}, s_k^i)], 
∀x_k A[\aSeqq{k-1}{s}{i}, x_k, f(\aSeqq{k-1}{s}{i}, x_k)], 
\Pi_i \fCenter \Lambda_i
$}
\LeftLabel{\nLu}
\RightLabel{\!.}
\UnaryInfC{$
∀x_k A[\aSeqq{k-1}{s}{i}, x_k, f(\aSeqq{k-1}{s}{i}, x_k)], 
\Pi_i \fCenter \Lambda_i
$}
\end{prooftree}
\def\proofSkipAmount{\vskip 0.8ex plus 0.8ex minus 0.4ex}
}
\noindent
Let us denote by \ARRi{s}{i} 
the sequence \aSeqq{k}{s}{i} of terms for $i = 1, \ldots, n$.
We may suppose without loss of generality that the sequence 
\(
\ARRi{s}{1}, \ldots, \ARRi{s}{n}
\)
of sequences of terms (possibly with repetitions) is arranged in such a way 
that for every $j < i$ the term $f(\ARRi{s}{j})$ 
is not a deep proper subterm of the term $f(\ARRi{s}{i})$.

Let $a_1, \ldots, a_{n}$ be the Henkin witnessing constants 
for the existential formulas
\(
∃y A[\ARRi{s}{1},y], \ldots, ∃y A[\ARRi{s}{n},y]
\):
\begin{align}
a_i ≡ \aHc{∃y A[\ARRi{s}{i},y]} \qquad (i = 1, \ldots, n)
.
\label{eq:sk:wf:thm:a}
\end{align}
We intend to eliminate the terms $f(\ARRi{s}{1}), \ldots, f(\ARRi{s}{n})$ from the proof \nP\
by (deeply) replacing them with the Henkin constants $a_1, \ldots, a_{n}$.
Let
\begin{align*}
\begin{aligned}
\ARRpi{s}{i} & ≡ 
\ARRi{s}{i} \aReplace{f(\ARRi{s}{1})}{a_1} \cdots \aReplace{f(\ARRi{s}{n})}{a_{n}}
\\
a'_i & ≡ 
f(\ARRi{s}{i}) \aReplace{f(\ARRi{s}{1})}{a_1} \cdots \aReplace{f(\ARRi{s}{n})}{a_{n}}
.
\end{aligned}
\quad (i = 1, \ldots, n)
\end{align*}
From the arrangement of terms in the sequence 
\(
\ARRi{s}{1}, \ldots, \ARRi{s}{n}
\)
we conclude that
\begin{align*}
\begin{aligned}
\ARRi{s}{i} & ≡ 
\ARRi{s}{i} \aReplace{f(\ARRi{s}{1})}{a_1} \cdots \aReplace{f(\ARRi{s}{i})}{a_{i}}
\\
a_i & ≡ 
f(\ARRi{s}{i}) \aReplace{f(\ARRi{s}{1})}{a_1} \cdots \aReplace{f(\ARRi{s}{i})}{a_{i}}
.
\end{aligned}
\quad (i = 1, \ldots, n)
\end{align*}
Consequently, $\ARRpi{s}{i}$ and $a'_i$ are obtained from $\ARRi{s}{i}$ and $a_i$ 
by deep replacement of only those terms $f(\ARRi{s}{j})$ by the corresponding $a_j$ 
for which the inequality $j > i$ holds:
\begin{align*}
\begin{aligned}
\ARRpi{s}{i} & ≡ 
\ARRi{s}{i} \aReplace{f(\ARRi{s}{i+1})}{a_{i+1}} \cdots \aReplace{f(\ARRi{s}{n})}{a_{n}}
\\
a'_i & ≡ 
a_i \aReplace{f(\ARRi{s}{i+1})}{a_{i+1}} \cdots \aReplace{f(\ARRi{s}{n})}{a_{n}}
.
\end{aligned}
\quad (i = 1, \ldots, n)
\end{align*}
(Hence $\ARRpi{s}{n} ≡ \ARRi{s}{n}$ and $a'_{n} ≡ a_{n}$.)
From this and \XEL{sk:wf:thm}{a} we immediately obtain that
\begin{align*}
a'_i ≡ \aHc{∃y A[\ARRpi{s}{i},y]} \qquad (i = 1, \ldots, n)
.
\end{align*}
The terms $a'_1, \ldots, a'_{n}$ are thus 
the Henkin witnessing constants for the existential formulas
\(
∃y A[\ARRpi{s}{1},y], \ldots, ∀y A[\ARRpi{s}{n},y]
\).

The following notation and terminology is used in the next elimination lemma.
By \nDr\ we denote the following composition of deep replacements:
\[
\nDr ≡ \aReplace{f(\ARRi{s}{1})}{a_1} \cdots \aReplace{f(\ARRi{s}{n})}{a_{n}}
.
\]
The above properties of deep replacements can be rewritten as follows:
\begin{align}
f(\ARRi{s}{i})\nDr ≡ a'_i ≡ \aHc{∃y A[\ARRi{s}{i}\nDr,y]} \qquad (i = 1, \ldots, n)
.
\label{eq:sk:wf:thm:theta:fs}
\end{align}
By \emph{special partial instances} of the one-step Skolemization \XEL{sk:wf:def}{sk} 
we mean formulas of the form ($1 \leq j \leq k$, $i = 1, \ldots, n$)
\begin{align}
\aSeqc[j]{k}{∀x} 
A[\aSeqq{j-1}{s}{i}, \aSeq[j]{k}{x}, f(\aSeqq{j-1}{s}{i}, \aSeq[j]{k}{x})]
.
\label{eq:sk:wf:thm:inst}
\end{align}
We use the symbol \nAwff\ to denote a set of formulas 
with special partial instances \XEL{sk:wf:thm}{inst} 
of \XEL{sk:wf:def}{sk} as its only elements.
Each formula \XEL{sk:wf:thm}{inst} 
is then transformed into the formula
\[
\aSeqc[j]{k}{∀x} ∃y A[\aSeqq{j-1}{s}{i}, \aSeq[j]{k}{x}, y]
.
\]
We denote by \nAwfe\ the image of \nAwff\ under this transformation.
\end{proof}

\paragraph{\textbf{Elimination lemma \textnormal{(Skolem witnessing function).}}}
{\myclaim

Let
\[
\nP[\rho] \rProofEq{l} \nAwff, \Pi \fCenter \Lambda
\]
be a proof free for the Skolem witnessing function $f$.
Let further 
\nAwff\ consists of some special partial instances of \XEL{sk:wf:def}{sk}
and let $m$ be the number of inferences in \nP[\rho] 
that are non-invariant under deep replacement of $f$-terms.
Then we can find a proof
\[
\aPp[\rho]{} \rProvableEq{l+m} \nAwfe\nDr, \Pi\nDr \fCenter \Lambda\nDr
\]
in $\nLG \cup \{f\}$
with the same cut-rank as \nP[\rho] 
invariant under deep replacement of $f$-terms.

}

\begin{proof}

By induction on the structure of the \nLKh-proof \nP[\rho].
We continue by case analysis of the last rule applied in \nP[\rho].
(As the reader will see, the proof transformation described below preserves 
cut-rank and yields a proof which is invariant under deep replacement of $f$-terms.
The verification of these easily checked facts is omitted here.)

Suppose first that the proof \nP[\rho]\ ends with an axiom of the form
\def\proofSkipAmount{\vskip 0.8ex plus 0.8ex minus 0.4ex}
\begin{prooftree}
\AxiomC{\phantom{axiom}}
\LeftLabel{\nAx}
\UnaryInfC{$\nAwff, B, \Pi_1 \fCenter \Lambda_1, B$}
\RightSubproofLabel{\nP[\rho]}
\end{prooftree}
\def\proofSkipAmount{\vskip 0.8ex plus 0.8ex minus 0.4ex}
with $B$ as principal formula.
Then it must be $l = 1$ and $m = 0$.
We built a new proof 
\def\proofSkipAmount{\vskip 0.8ex plus 0.8ex minus 0.4ex}
\begin{prooftree}
\AxiomC{\phantom{axiom}}
\LeftLabel{\nAx}
\UnaryInfC{$\nAwfe\nDr, B\nDr, \Pi_1\nDr \fCenter \Lambda_1\nDr, B\nDr$}
\RightSubproofLabel{\nPp[\rho]}
\end{prooftree}
\def\proofSkipAmount{\vskip 0.8ex plus 0.8ex minus 0.4ex}
with $B\nDr$ as principal formula.
We clearly have
\(
\nPp[\rho] \rProofEq{1+0} \nAwfe\nDr, B\nDr, \Pi_1\nDr \fCenter \Lambda_1\nDr, B\nDr
\).

Suppose now that the proof \nP[\rho]\ ends with an \nLu-inference of the form
\def\proofSkipAmount{\vskip 0.8ex plus 0.8ex minus 0.4ex}
\begin{prooftree}
\AxiomC{}
\RightLabel{\aP[\rho]{1}}
\branchDeduce
\DeduceC{$
\nAwff, A[\aSeqnq{k-1}{s}{i}, \aT{s}{k}{\nIs}, f(\aSeqnq{k-1}{s}{i}, \aT{s}{k}{\nIs})],
\Pi \fCenter \Lambda
$}
\LeftLabel{\nLu}
\RightLabel{,}
\UnaryInfC{$\nAwff, \Pi \fCenter \Lambda$}
\RightSubproofLabel{\nP[\rho]}
\end{prooftree}
\def\proofSkipAmount{\vskip 0.8ex plus 0.8ex minus 0.4ex}
where the principal formula of the \nLu-rule is
\[
∀x_k A[\aSeqnq{k-1}{s}{i}, x_k, f(\aSeqnq{k-1}{s}{i}, x_k)] \in \nAwff
.
\]
We have 
\(
\aSeqnq{k}{s}{\nIs} ≡ \ARRi{s}{i} 
\)
for some $i = 1, \ldots, n$.
We apply the inductive hypothesis to the subderivation of \nP[\rho]:
\[
\aP[\rho]{1} \rProofEq{l-1} 
\nAwff, A[\aSeqnq{k-1}{s}{i}, \aT{s}{k}{\nIs}, f(\aSeqnq{k-1}{s}{i}, \aT{s}{k}{\nIs})],
\Pi \fCenter \Lambda
\]
and find a proof \aPp[\rho]{1} such that
\[
\aPp[\rho]{1} \rProofEq{l-1+m-1}
\nAwfe\nDr,
A[\aSeqnqdr{k-1}{s}{i}{\nDr}, \aT{s}{k}{\nIs}\nDr, 
  f(\aSeqnq{k-1}{s}{i}, \aT{s}{k}{\nIs})\nDr],
\Pi\nDr \fCenter \Lambda\nDr
.
\]
From \XEL{sk:wf:thm}{theta:fs} we can see that the term
\[
f(\aSeqnq{k-1}{s}{i}, \aT{s}{k}{\nIs})\nDr ≡ 
\aHc{∃y A[\aSeqnqdr{k-1}{s}{i}{\nDr}, \aT{s}{k}{\nIs}\nDr, y]} ≡ a'
\]
is, in fact, the Henkin witnessing constant for
\(
∃y A[\aSeqnqdr{k-1}{s}{i}{\nDr}, \aT{s}{k}{\nIs}\nDr, y]
\).
Hence
\begin{align*}
\aPp[\rho]{1} \rProofEq{l+m-2}
\nAwfe\nDr, A[\aSeqnqdr{k-1}{s}{i}{\nDr}, \aT{s}{k}{\nIs}\nDr, a'],
\Pi\nDr \fCenter \Lambda\nDr
.
\end{align*}
We now construct a new proof
\def\proofSkipAmount{\vskip 0.8ex plus 0.8ex minus 0.4ex}
\begin{prooftree}
\AxiomC{}
\RightLabel{\aPp[\rho]{1}}
\branchDeduce
\DeduceC{$
\nAwfe\nDr, A[\aSeqnqdr{k-1}{s}{i}{\nDr}, \aT{s}{k}{\nIs}\nDr, a'],
\Pi\nDr \fCenter \Lambda\nDr
$}
\LeftLabel{\nLe}
\UnaryInfC{$
\nAwfe\nDr, ∃y A[\aSeqnqdr{k-1}{s}{i}{\nDr}, \aT{s}{k}{\nIs}\nDr, y],
\Pi\nDr \fCenter \Lambda\nDr
$}
\LeftLabel{\nLu}
\RightLabel{.}
\UnaryInfC{$\nAwfe\nDr, \Pi\nDr \fCenter \Lambda\nDr$}
\RightSubproofLabel{\nPp[\rho]}
\end{prooftree}
\def\proofSkipAmount{\vskip 0.8ex plus 0.8ex minus 0.4ex}
The last step is an \nLu-inference with the principal formula
\[
∀x_k ∃y A[\aSeqnqdr{k-1}{s}{i}{\nDr}, x_k, y] \in \nAwfe\nDr
.
\]
We obviously have
\(
\nPp[\rho] \rProofEq{l+m} \nAwfe\nDr, \Pi\nDr \fCenter \Lambda\nDr
\).

Suppose that the proof \nP[\rho]\ ends with an \nLu-inference of the form
($1 \leq i < k$)
{\relsize{-0.5}
\def\proofSkipAmount{\vskip -1.0ex plus 0.8ex minus 0.4ex}
\begin{prooftree}
\hskip -2.5em
\AxiomC{}
\RightLabel{\aP[\rho]{1}}
\branchDeduce
\DeduceC{$
\aSeqcdd[\nIq+1]{k}{∀x} 
A[\aSeqnqdd{\nIq-1}{s}{i}, \aT{s}{\nIq}{\nIs}, \aSeqdd[\nIq+1]{k}{x}, 
  f(\aSeqnqdd{\nIq-1}{s}{i}, \aT{s}{\nIq}{\nIs}, \aSeqdd[\nIq+1]{k}{x})],
\nAwff, \Pi \fCenter \Lambda
$}
\LeftLabel{\nLu}
\RightLabel{,}
\UnaryInfC{$\nAwff, \Pi \fCenter \Lambda$}
\RightSubproofLabel{\nP[\rho]}
\end{prooftree}
\def\proofSkipAmount{\vskip 0.8ex plus 0.8ex minus 0.4ex}
}
\noindent
where the principal formula of the \nLu-rule is
\[
∀x_\nIq \aSeqc[\nIq+1]{k}{∀x} 
A[\aSeqnq{\nIq-1}{s}{i}, x_\nIq, \aSeq[\nIq+1]{k}{x}, 
  f(\aSeqnq{\nIq-1}{s}{i}, x_\nIq, \aSeq[\nIq+1]{k}{x})]
\in \nAwff
.
\]
The sequence of terms \aSeqnq{\nIq}{s}{\nIs} is an initial segment of some sequence from
\(
\ARRi{s}{1}, \ldots, \ARRi{s}{n}
\).
We apply the inductive hypothesis to the subderivation of \nP[\rho]:
{\relsize{-0.5}
\begin{align*}
\aP[\rho]{1} \rProofEq{l-1} 
\aSeqc[\nIq+1]{k}{∀x} 
A[\aSeqnq{\nIq-1}{s}{i}, \aT{s}{\nIq}{\nIs}, \aSeq[\nIq+1]{k}{x}, 
  f(\aSeqnq{\nIq-1}{s}{i}, \aT{s}{\nIq}{\nIs}, \aSeq[\nIq+1]{k}{x})],
\nAwff, \Pi \fCenter \Lambda
\end{align*}
}%
and find a proof \aPp[\rho]{1} such that
\[
\aPp[\rho]{1} \rProofEq{l-1+m}
\aSeqc[\nIq+1]{k}{∀x} ∃y
A[\aSeqnqdr{\nIq-1}{s}{i}{\nDr}, \aT{s}{\nIq}{\nIs}\nDr, \aSeq[\nIq+1]{k}{x}, y],
\nAwfe\nDr, \Pi\nDr \fCenter \Lambda\nDr
.
\]
We use an \nLu-inference with the principal formula
\[
∀x_\nIq \aSeqc[\nIq+1]{k}{∀x} ∃y
A[\aSeqnqdr{\nIq-1}{s}{i}{\nDr}, x_\nIq, \aSeq[\nIq+1]{k}{x}, y]
\in \nAwfe\nDr
\]
and construct a new proof
\def\proofSkipAmount{\vskip -2.5ex plus 0.8ex minus 0.4ex}
\begin{prooftree}
\hskip -3.0em
\AxiomC{}
\RightLabel{\aPp[\rho]{1}}
\branchDeduce
\DeduceC{$
\aSeqc[\nIq+1]{k}{∀x} ∃y
A[\aSeqnqdr{\nIq-1}{s}{i}{\nDr}, \aT{s}{\nIq}{\nIs}\nDr, \aSeq[\nIq+1]{k}{x}, y],
\nAwfe\nDr, \Pi\nDr \fCenter \Lambda\nDr
$}
\LeftLabel{\nLu}
\RightLabel{.}
\UnaryInfC{$\nAwfe\nDr, \Pi\nDr \fCenter \Lambda\nDr$}
\RightSubproofLabel{\nPp[\rho]}
\end{prooftree}
\def\proofSkipAmount{\vskip 0.8ex plus 0.8ex minus 0.4ex}
We thus have
\(
\nPp[\rho] \rProofEq{l+m} \nAwfe\nDr, \Pi\nDr \fCenter \Lambda\nDr
\).

Suppose that the proof \nP[\rho]\ ends with an \nRu-inference of the form
\def\proofSkipAmount{\vskip -1.0ex plus 0.8ex minus 0.4ex}
\begin{prooftree}
\AxiomC{}
\RightLabel{\aP[\rho]{1}}
\branchDeduce
\DeduceC{$\nAwff, \Pi \fCenter \Lambda_1, \aSubst{B}{z}{\aHc{∀z B}}$}
\LeftLabel{\nRu}
\RightLabel{.}
\UnaryInfC{$\nAwff, \Pi \fCenter \Lambda_1, ∀z B$}
\RightSubproofLabel{\nP[\rho]}
\end{prooftree}
\def\proofSkipAmount{\vskip 0.8ex plus 0.8ex minus 0.4ex}
By assumption,
the proof \aP[\rho]{} is free for the Skolem witnessing function $f$.
Hence the \nRu-inference is invariant under deep replacement of $f$-terms
and thus
\[
\aSubst{B}{z}{\aHc{∀z B}}\nDr ≡ \aSubst{B\nDr}{z}{\aHc{∀z B\nDr}}
.
\]
We apply the inductive hypothesis to the subderivation
\(
\aP[\rho]{1} \rProofEq{l-1} 
\nAwff, \Pi \fCenter \Lambda_1, \aSubst{B}{z}{\aHc{∀z B}}
\)
of \aP[\rho]{} and find a proof \aPp[\rho]{1} such that
\[
\aPp[\rho]{1} \rProofEq{l-1+m}
\nAwfe\nDr, \Pi\nDr \fCenter \Lambda_1\nDr, \aSubst{B\nDr}{z}{\aHc{∀z B\nDr}}
.
\]
We now use a similar \nRu-inference and construct a proof
\def\proofSkipAmount{\vskip -1.0ex plus 0.8ex minus 0.4ex}
\begin{prooftree}
\AxiomC{}
\RightLabel{\aPp[\rho]{1}}
\branchDeduce
\DeduceC{$
\nAwfe\nDr, \Pi\nDr \fCenter \Lambda_1\nDr, \aSubst{B\nDr}{z}{\aHc{∀z B\nDr}}
$}
\LeftLabel{\nRu}
\RightLabel{.}
\UnaryInfC{$
\nAwfe\nDr, \Pi\nDr \fCenter \Lambda_1\nDr, ∀z B\nDr
$}
\RightSubproofLabel{\nPp[\rho]}
\end{prooftree}
\def\proofSkipAmount{\vskip 0.8ex plus 0.8ex minus 0.4ex}
We clearly have
\(
\nPp[\rho] \rProofEq{l+m} 
\nAwfe\nDr, \Pi\nDr \fCenter \Lambda_1\nDr, ∀z B\nDr
\).

Suppose that the proof \nP[\rho]\ ends with a cut inference of the form
\def\proofSkipAmount{\vskip -1.0ex plus 0.8ex minus 0.4ex}
\begin{prooftree}
\AxiomC{}
\RightLabel{\aP[\rho]{1}}
\branchDeduce
\DeduceC{$\nAwff, \Pi \fCenter \Lambda, B$}
\AxiomC{}
\RightLabel{\aP[\rho]{2}}
\branchDeduce
\DeduceC{$B, \nAwff, \Pi \fCenter \Lambda$}
\LeftLabel{\nRa}
\RightLabel{,}
\BinaryInfC{$\nAwff, \Pi \fCenter \Lambda$}
\RightSubproofLabel{\nP[\rho]}
\end{prooftree}
\def\proofSkipAmount{\vskip 0.8ex plus 0.8ex minus 0.4ex}
where
\(
\aP[\rho]{1} \rProofEq{l_1} \nAwff, \Pi \fCenter \Lambda, B
\)
and
\(
\aP[\rho]{2} \rProofEq{l_2} \nAwff, B, \Pi \fCenter \Lambda
\)
for some numbers $l_1, l_2$ such that $l_1 + l_2 + 1 = l$.
We may assume w.l.o.g. that the subproof \aP[\rho]{1} contains 
the first $m_1$ non-invariant \nLu-rules and
the subproof \aP[\rho]{2} the remaining $m_2 = m-m_1$.
We apply the inductive hypothesis to both subderivations and obtain
\begin{align*}
\aPp[\rho]{1} \rProofEq{l_1+m_1} 
\nAwfe\nDr, \Pi\nDr \fCenter \Lambda\nDr, B\nDr
\quad \text{and} \quad
\aPp[\rho]{2} \rProofEq{l_2+m_2} 
\nAwfe\nDr, B\nDr, \Pi\nDr \fCenter \Lambda\nDr
.
\end{align*}
We use a similar cut rule and built a new proof
\def\proofSkipAmount{\vskip -1.0ex plus 0.8ex minus 0.4ex}
\begin{prooftree}
\AxiomC{}
\RightLabel{\aPp[\rho]{1}}
\branchDeduce
\DeduceC{$\nAwfe\nDr, \Pi\nDr \fCenter \Lambda\nDr, B\nDr$}
\AxiomC{}
\RightLabel{\aPp[\rho]{2}}
\branchDeduce
\DeduceC{$B\nDr, \nAwfe\nDr, \Pi\nDr \fCenter \Lambda\nDr$}
\LeftLabel{\makebox[0cm][r]{\nRa}}
\RightLabel{.}
\BinaryInfC{$\nAwfe\nDr, \Pi\nDr \fCenter \Lambda\nDr$}
\RightSubproofLabel{\nPp[\rho]}
\end{prooftree}
\def\proofSkipAmount{\vskip 0.8ex plus 0.8ex minus 0.4ex}
From the identity
$l_1 + m_1  + l_2 + m_2 + 1 = l_1 + l_2 + 1 + m_1 + m_2 = l + m$
we conclude that
\(
\nPp[\rho] \rProofEq{l+m}
\nAwfe\nDr, \Pi\nDr \fCenter \Lambda\nDr 
\).

The remaining cases are proved similarly.
This proves the elimination lemma.
\end{proof}

\begin{proof}[Proof of the elimination theorem continued]

By the elimination lemma, where \nAwff\ is a singleton multiset
\(
\aSeqc{k}{∀x} A[\aSeq{k}{x},f(\aSeq{k}{x})]
\),
we obtain a proof in $\nLG \cup \{f\}$
\[
\aPpp{} \rProofEq{l+n}
\aSeqc{k}{∀x} ∃y A[\aSeq{k}{x},y], \Gamma \fCenter \Delta
\]
with the same cut-rank as \aP{} 
invariant under deep replacement of $f$-terms.
We eliminate the remaining $f$-terms by deeply replacing them 
by some terms from the Henkin expansion of \nLG.
Therefore, by the replacement lemma~\ref{par:repl:thm},
we can find a proof 
\begin{align*}
\aPp{} \rProofEq{l+n} 
\aSeqc{k}{∀x} ∃y A[\aSeq{k}{x},y], \Gamma \fCenter \Delta
\end{align*}
\emph{in \nLG}
with the same cut-rank as \aP{}.
This proves the elimination theorem \ref{par:sk:wf:thm}.
\end{proof}

\PAR{Corollary}
\label{par:sk:wf:cor}
{\myclaim

Let $B$ be a \nLG-formula.
If the formula $∀\ARR{x} A[\ARR{x},f(\ARR{x})] → B$ has an \nLKh-proof 
free for the Skolem witnessing function $f$ 
with length $l$ and cut-rank $≤ r$,
then the formula $∀\ARR{x} ∃y A → B$ has an \nLKh-proof in \nLG\ 
with length $\leq 2l$ and cut-rank $≤ r$.

}

\begin{proof}
Let
\(
\aP{} \rProofEq{l} ∀x A[\ARR{x},f(\ARR{x})] → B
\)
be a proof free for $f$ with cut-rank $≤ r$.
By the inversion lemma \ref{par:pr:inv},
since inversion preserves the property of being free 
for a witnessing Skolem function, 
we can find a proof
\(
\aP{1} \rProofLe{l} ∀x A[\ARR{x},f(\ARR{x})] \fCenter B
\)
free for $f$ with cut-rank $≤ r$.
By applying the elimination theorem \ref{par:sk:wf:thm} and 
$l+n < l+l = 2l$, we have
\(
\rProofLt{2l} ∀\ARR{x} ∃y A \fCenter B
\).
By applying one \nRi-inference, we obtain 
\(
\rProofLe{2l} ∀\ARR{x} ∃y A → B
\).
\end{proof}

\PAR{Problem 22}
\label{par:sk1:problem22}

The original Problem~22 by P.~Pudlák listed in~\cite{qq:CloteKrajicek:93}
is the following question (using our notation and terminology):
\begin{quote}
Assume that $∀\ARRr{x}∃yA[\ARRr{x},y]$ is provable in predicate logic.
Introduce a new function symbol $f$ and 
an axiom that states $∀\ARRr{x}A[\ARRr{x},f(\ARRr{x})]$.
Does there exist semiformula $A$ such that the extended system 
gives a superexponential speed-up over predicate calculus
with respect to number of symbols in proofs?
\end{quote}
We give a negative answer to the problem 
for a subclass of \nLKh-proofs consisting of \nLKh-derivations 
that are free for the Skolem witnessing function $f$ which belongs 
to the one-step Skolemization $∀\ARRr{x}A[\ARRr{x},f(\ARRr{x})]$ 
of the formula $∀\ARRr{x}∃yA[\ARRr{x},y]$.
This fact is a straightforward consequence of the next theorem.

The modified subtraction function $x \bMinus y$ on natural numbers is defined by
\begin{align*}
x \bMinus y = z ↔ x ≥ y ∧ z+y = x ∨ x ≤ y ∧ z = 0
.
\end{align*}
Recall also that
\(
\aHtP{\nP} ≤ l < 2^{\aHtP{\nP}}
\)
if the length of the \nLKh-proof \nP\ is $l$.

\PAR{Theorem}
\label{par:sk1:problem22:thm}
{\myclaim

Let
\(
\aP{A} \rProof ∀\ARR{x} ∃y A[\ARR{x},y]
\)
be an \nLKh-proof and $f$ a new function symbol.
Then an elementary function \aE{l}{r} exists, which depends only on \aP{A},
such that the following holds.
If
\(
\nP \rProofEq{l} ∀\ARR{x} A[\ARR{x},f(\ARR{x})] → B
\)
is an \nLKh-proof with cut-rank $≤ r$
free for the Skolem witnessing function $f$ 
that belongs to the one-step Skolemization 
$∀\ARR{x}A[\ARR{x},f(\ARR{x})]$ of $∀\ARR{x}∃yA[\ARR{x},y]$,
then there is an \nLKh-proof
\(
\nP[\rho] \rProofLe{\aE{l}{r}} B
\)
with cut-rank $≤ r$.

}

\begin{proof}
Set
\(
d_A = \aDpF{∀\ARR{x}∃yA}
\),
\(
r_A = \aCrP{\aP{A}}
\)
and
\(
h_A = \aHtP{\aP{A}}
\).
By Cor.~\ref{par:sk:wf:cor} and 
the inversion lemma~\ref{par:pr:inv} applied to \aP{},
we can find a proof
\(
\nPp \rProofLe{2l} ∀\ARR{x} ∃y A \fCenter B
\)
with cut-rank $≤ r$.
We build a new proof using one weakening 
that ends in a cut on the formula $∀\ARR{x}∃yA$:
\def\proofSkipAmount{\vskip 0.8ex plus 0.8ex minus 0.4ex}
\begin{prooftree}
\AxiomC{}
\RightLabel{\aWkP{\aP{A}}{}{B, ∀\ARR{x}∃yA}}
\branchDeduce
\DeduceC{$\fCenter B, ∀\ARR{x}∃yA$}
\AxiomC{}
\RightLabel{\nPp}
\branchDeduce
\DeduceC{$∀\ARR{x}∃yA \fCenter B$}
\LeftLabel{\nCut}
\RightLabel{.}
\BinaryInfC{$\fCenter B$}
\RightSubproofLabel{\aPp[\rho]{}}
\end{prooftree}
\def\proofSkipAmount{\vskip 0.8ex plus 0.8ex minus 0.4ex}
The proof 
\(
\aPp[\rho]{} \rProof B
\)
has the cut-rank 
\begin{align*}
\aCrP{\aPp[\rho]{}} & =
\max\left(d_A,\aCrP{\aWkP{\aP{A}}{}{B, ∀\ARR{x}∃yA}},\aCrP{\aPp{}}\right) 
 \mEQx{L \ref{par:pr:weak}}{=}
\max\left(d_A,r_A,\aCrP{\aPp{}}\right) 
\\ & ≤
\max\left(d_A,r_A,r\right) =
\max\left(d_A,r_A\right) \bMinus r + r
\end{align*}
and the height
\begin{align*}
\aHtP{\aPp[\rho]{}} & =
\max\left(\aHtP{\aWkP{\aP{A}}{}{B, ∀\ARR{x}∃yA}},\aHtP{\aPp{}}\right)+1 
 \mEQx{L \ref{par:pr:weak}}{=}
\max\left(h_A,\aHtP{\aPp{}}\right)+1 ≤
\max\left(h_A,2l\right)+1
.
\end{align*}
By the $n$-fold application of the reduction lemma~\ref{par:cut:red:III}, where
\(
n = \max\left(d_A,r_A\right) \bMinus r
\),
we can find a proof
\(
\aP[\rho]{} \rProof B
\)
with cut-rank $≤ r$ and with height
\begin{align*}
\aHtP{\aP[\rho]{}} ≤
2^{\aHtP{\aPp[\rho]{}}}_{2\max\left(d_A,r_A\right)\bMinus2r} ≤
2^{\max\left(h_A,2l\right)+1}_{2\max\left(d_A,r_A\right)\bMinus2r}
.
\end{align*}
The length of the proof \aP[\rho]{} 
is thus strictly bounded by the number
\begin{align*}
2^{\aHtP{\aP[\rho]{}}} ≤
2^{\max\left(h_A,2l\right)+1}_{2\max\left(d_A,r_A\right)\bMinus2r+1} ≤
2^{\max\left(h_A,2l\right)+1}_{2\max\left(d_A,r_A\right)+1}
.
\end{align*}
It suffices to define \nE\ by
\begin{align*}
\aE{l}{r} = 
2^{\max\left(h_A,2l\right)+1}_{2\max\left(d_A,r_A\right)+1} 
\end{align*}
as an elementary function.
This concludes the proof since obviously
\(
\nP[\rho] \rProofLe{\aE{l}{r}} B
\).
\end{proof}

\PAR{Skolem counterexample functions}
\label{par:sk:cf:def}

Let \nLG\ be a first-order language and $g$ a new $k$-ary function symbol not \nLG.
Let further
\begin{align}
\aSeqc{k}{∃x}∀yA[\aSeq{k}{x},y]
\label{eq:sk:cf:def:f}
\end{align}
be an \nLG-formula. 
A \emph{one-step counterexample Skolemization} of 
\XEL{sk:cf:def}{f} 
is a formula
\begin{align}
\aSeqc{k}{∃x} A[\aSeq{k}{x}, g(\aSeq{k}{x})]
\label{eq:sk:cf:def:sk}
.
\end{align}
The function $g$ is called the \emph{Skolem counterexample function} of the skolemization.

We say that an \nLKh-proof \nP\ is 
\emph{free for the Skolem counterexample function $g$} 
(\emph{free for $g$} for short)
if the only rules in the proof \nP\
that are non-invariant under deep replacement of $g$-terms
are \nRe-inferences of the form
{\relsize{-0.5}
\def\proofSkipAmount{\vskip -0.1ex plus 0.8ex minus 0.4ex}
\begin{prooftree}
\AxiomC{$
\Pi \fCenter \Lambda,
∃x_k A[\aSeqnq{k-1}{s}{i}, x_k, g(\aSeqnq{k-1}{s}{i}, x_k)],
A[\aSeqnq{k-1}{s}{i}, s_k, g(\aSeqnq{k-1}{s}{i}, s_k)]
$}
\LeftLabel{\nRe}
\RightLabel{.}
\UnaryInfC{$
\Pi \fCenter \Lambda,
∃x_k A[\aSeqnq{k-1}{s}{i}, x_k, g(\aSeqnq{k-1}{s}{i}, x_k)]
$}
\end{prooftree}
\def\proofSkipAmount{\vskip 0.8ex plus 0.8ex minus 0.4ex}
}
\noindent
By the subformula property~\ref{par:cut:sub} cut-free proofs are trivially free for $g$.
The property of an \nLKh-proof being free for $g$ is downwards hereditary property:
subderivations of such a proof are themselves free for $g$.

In the following theorem we assume that 
the first-order language \nLG, 
the Skolem counterexample function $g$ and
the skolemization formula
are the same as above.
Note that we are working in the Henkin expansion 
of the first-order language $\nLG \cup \{g\}$.

\PAR{Elimination theorem \textnormal{(Skolem counterexample function)}}
\label{par:sk:cf:thm}
{\myclaim

Let $\Gamma$ and $\Delta$ be finite multisets of \nLG-formulas.
Let further
\[
\nP \rProofEq{l} \Gamma \fCenter \Delta, \aSeqc{k}{∃x} A[\aSeq{k}{x},g(\aSeq{k}{x})]
\]
be an \nLKh-proof free for the Skolem counterexample function $g$.
Let finally $n$ be the number of inferences in \aP{} 
which are non-invariant under deep replacement of $g$-terms.
Then there is an \nLKh-proof \nPp\ in \nLG\ 
with the same cut-rank as \nP\ such that
\[
\nPp \rProofEq{l+n} \Gamma \fCenter \Delta, \aSeqc{k}{∃x} ∀y A[\aSeq{k}{x},y)]
.
\]

}

\begin{proof}
By ``dualizing'' the arguments which are used 
in the proof of the elimination theorem~\ref{par:sk:wf:thm}.
\end{proof}

\PAR{Corollary}
\label{par:sk:cf:cor}
{\myclaim

Let
\begin{align}
∃\ARRi{x}{1} ∀y_{1} ∃\ARRi{x}{2} ∀y_{2} \ldots ∃\ARRi{x}{n} ∀y_{n}
A[\ARRi{x}{1}, y_{1}, \ARRi{x}{2}, y_{2}, \ldots, \ARRi{x}{n},y_{n}]
\label{eq:sk:cf:cor:f}
\end{align}
be an \nLG-formula.
Let further $g_1, g_2, \ldots, g_n$ be pairwise different function symbols
of appropriate arities not in \nLG\ such that
each function $g_i$ for $i=1,\ldots,n$ is the Skolem counterexample function 
belonging to the one-step Skolemization 
{\relsize{-0.5}
\begin{align*}
∃\ARRi{x}{1}
\ldots ∃\ARRi{x}{i}
∃\ARRi{x}{i+1} ∀y_{i+1} \ldots ∃\ARRi{x}{n} ∀y_n 
A[\ARRi{x}{1}, g_1(\ARRi{x}{1}), \ldots,
  \ARRi{x}{i}, g_i(\ARRi{x}{1}, \ldots, \ARRi{x}{i}), 
  \ARRi{x}{i+1}, y_{i+1}, \ldots, \ARRi{x}{n}, y_{n}
  ]
\end{align*}
}%
of the formula 
{\relsize{-0.5}
\begin{align*}
∃\ARRi{x}{1}
\ldots ∃\ARRi{x}{i} ∀y_i
∃\ARRi{x}{i+1} ∀y_{i+1} \ldots ∃\ARRi{x}{n} ∀y_n 
A[\ARRi{x}{1}, g_1(\ARRi{x}{1}), \ldots,
  \ARRi{x}{i}, y_i, 
  \ARRi{x}{i+1}, y_{i+1}, \ldots, \ARRi{x}{n}, y_{n}
  ]
.
\end{align*}
}%
If the formula 
\begin{align}
∃\ARRi{x}{1} ∃\ARRi{x}{2} \ldots ∃\ARRi{x}{n} 
A[\ARRi{x}{1}, g_1(\ARRi{x}{1}),
  \ARRi{x}{2}, g_2(\ARRi{x}{1}, \ARRi{x}{2}), \ldots,
  \ARRi{x}{n}, g_n(\ARRi{x}{1}, \ARRi{x}{2}, \ldots, \ARRi{x}{n})]
\label{eq:sk:cf:cor:sk}
\end{align}
has an \nLKh-proof free for the Skolem counterexample functions \aSeq{n}{g}
with length $l$ and with cut-rank $≤ r$,
then the formula \XEL{sk:cf:cor}{f}
has an \nLKh-proof in \nLG\ with length $< 2l$ and with cut-rank $≤ r$.

}

\begin{proof}
By $n$-fold application of the elimination theorem \ref{par:sk:cf:thm}
starting with $g_1$ and ending with $g_n$.
Elimination of one Skolem function preserves 
the number of rules non-invariant under deep replacement of Skolem terms
for the remaining functions.
\end{proof}

\paragraph{Remark.} 

We do not require for the semiformula $A$ in \XEL{sk:cf:cor}{f} to be quantifier-free.
Hence the quantifier formula \XEL{sk:cf:cor}{f} might not to be in prenex normal form
(cf. \cite[Prop.~1]{DBLP:journals/jsyml/BaazHW12}.)


\section{Conclusion}
\label{sc:co}

We have given a partial positive solution to the problem 
of efficient elimination of Skolem functions from proofs 
in first-order logic without equality:
elimination of a single Skolem function from \nLKh-proofs with cuts
that are free for the Skolem function 
increases the length of such proofs only linearly.
The question is whether we can lift the restrictions 
imposed upon cut inferences  in \nLKh-derivations and 
allow the use of more general lemmas that 
contain Skolem semiterms in the scope of quantifiers.
The preliminary results in this direction are looking promising.

Our second goal has been to introduce 
a new variant of Gentzen's sequent calculus that we call \nLKh. 
Let us summarize some of its specific distinguished features:
\begin{itemize}
\item
Local soundness:
\nLKh-inferences preserve validity in Henkin structures.
\item
The characterization problem:
\nLKh-provable formulas, 
even with parameters, 
are exactly those that are valid in every Henkin structure.
\item
Strong locality property:
validity of each \nLKh-inference step depends only
on the form of its active formulas.
\item
Deep replacement of terms:
this self-correcting local transformation turns each invariant \nLKh-rule 
into another valid inference of the same kind.
\item
Effective cut elimination 
via a non-Gentzen-style algorithm without resorting to regularization.
\end{itemize}
Most of these properties has been already established for its kin, a tableaux-based calculus from 
Komara and Voda \cite{qq:KomaraVoda:95:TABLEAUX95}.
As the calculi \nLKh\ and 
\nLKpp\ of Aguilera and Baaz~\cite{DBLP:journals/jsyml/AguileraB19}
simulate each other linearly,
proof-theoretic results acquired for one system 
can be easily applied to the other, and vice versa.
For instance, we have:
\begin{itemize}
\item
Nonelementary speed-up:
\nLKh\ yields non-elementarily shorter cut-free proofs than the standard \nLK-proofs.
\item
Effective cut elimination for \nLKpp-proofs
via a non-Gentzen-style algorithm without resorting to regularization.
\end{itemize}
The first property follows directly 
from \cite[Cor. 2.7]{DBLP:journals/jsyml/AguileraB19},
originally proved for \nLKpp,
while the second one from the corresponding property for \nLKh.

The listed properties of \nLKh\ usually make structural proof analysis easier,
sometimes even with a better proof-complexity result.
The witness is the efficient elimination of Skolem functions.
Is it worth to study such a system? Only time will tell.



\newcommand{\noopsort}[1]{} \newcommand{\singletter}[1]{#1}
\providecommand{\href}[2]{#2}\begingroup\raggedright\endgroup


\end{document}